\documentclass[preprint,5p,times,twocolumn]{elsarticle}

\usepackage{amsmath,amsthm,amssymb,mathtools,mathrsfs} 
\allowdisplaybreaks[4]
\usepackage{tikz}
\usetikzlibrary{cd,patterns}
\usepackage{graphicx}
\usepackage{placeins}
\usepackage{cuted}
\usepackage{subcaption}
\usepackage{flushend}

\usepackage[textsize=tiny,textwidth=1.3cm]{todonotes}

\usepackage[hidelinks,
            colorlinks=true,
            linkcolor=blue,
            filecolor=magenta,      
            urlcolor=cyan,
            citecolor = blue
            ]{hyperref}

\newtheorem{theorem}{Theorem}[section]

\newtheorem{lemma}{Lemma}[section]
\newtheorem{proposition}{Proposition}[section]

\theoremstyle{definition}

\newtheorem{remark}{Remark}

\usepackage{enumitem}
\setenumerate[1]{itemsep=0pt,partopsep=0pt,parsep=\parskip,topsep=2.5pt}
\setenumerate[2]{itemsep=0pt,partopsep=0pt,parsep=\parskip,topsep=2.5pt}
\setenumerate[3]{itemsep=0pt,partopsep=0pt,parsep=\parskip,topsep=2.5pt}

\setitemize[1]{itemsep=0pt,partopsep=0pt,parsep=\parskip,topsep=2.5pt}
\setitemize[2]{itemsep=0pt,partopsep=0pt,parsep=\parskip,topsep=2.5pt}
\setitemize[3]{itemsep=0pt,partopsep=0pt,parsep=\parskip,topsep=2.5pt}

\setdescription[1]{itemsep=0pt,partopsep=0pt,parsep=\parskip,topsep=2.5pt}
\setdescription[2]{itemsep=0pt,partopsep=0pt,parsep=\parskip,topsep=2.5pt}
\setdescription[3]{itemsep=0pt,partopsep=0pt,parsep=\parskip,topsep=2.5pt}

\usepackage{fancyhdr} 
\usepackage{blindtext}

\let\emph\relax 
\DeclareTextFontCommand{\emph}{\bfseries\em}

\usepackage{etoolbox}
\makeatletter
\def\ps@pprintTitle{%
  \let\@oddhead\@empty
  \let\@evenhead\@empty
  \def\@oddfoot{}%
  \let\@evenfoot\@oddfoot}
\makeatother

\begin{document}

\begin{frontmatter}
    
\title{Multi-Target Tracking via Mean Field Feedback Control}
\author[1]{Changming Lei}
\ead{changming_lei@bit.edu.cn}

\author[1]{Junhao Dai}
\ead{a1446162133@gmail.com}

\author[1]{Donghua Shi\corref{cor1}}
\cortext[cor1]{Corresponding author}
\ead{dshi@bit.edu.cn}

\affiliation[1]{
  organization={School of Mathematics and Statistics, Beijing Institute of Technology},
  city={Beijing},
  postcode={100081},
  country={China}
}

\begin{abstract}
Cooperative tracking of multiple moving targets requires a swarm to continually adjust its deployment. Existing methods often rely on repeated agent–target reassignment or updates to position-dependent partitions, which can complicate coordination as swarms grow. This motivates us to seek an approach with a lighter coordination burden.
Specifically, we formulate multi-target tracking as a time-varying density-tracking problem, with the targets' tracking and sensing demands encoded in a task-demand density. From an Onsager-type variational formulation, a closed-form mean-field feedback law is derived to make the swarm density track this evolving density, and the resulting closed-loop system is shown to be well posed. Under the stated assumptions, a Kullback--Leibler (KL) divergence analysis establishes exponential KL convergence under ideal feedback and conditional robustness to feedback mismatch under estimated.
Numerical experiments demonstrate effective tracking of multiple moving targets and show smaller density errors at comparable applied control energy compared with an existing density-tracking feedback method.
\end{abstract}
\begin{keyword}
mean-field feedback control \sep density tracking \sep Onsager variational principle \sep Kullback--Leibler divergence
\end{keyword}
\end{frontmatter}

\section{Introduction}

Multi-agent systems (MASs) are increasingly used in persistent sensing, environmental monitoring, search and rescue, and precision agriculture~\cite{Brambilla2013SwarmReview,Dorigo2021SwarmRobotics,Elamvazhuthi2020Survey}. Among the tasks arising in these settings, multi-target tracking is a fundamental task of considerable practical importance in robotics~\cite{Robin2016MultiRobotTracking}. 

In this study, we consider a multi-target tracking problem in which a swarm is required to follow multiple moving targets, allocate sensing resources among them, and provide the desired spatial coverage for each target. At the individual-agent level, previous studies have addressed this problem through two main approaches. Explicit task-allocation methods establish robot--target assignments \cite{Gerkey2004TaskAllocation,Korsah2013TaskAllocation}. Coverage-control methods, by contrast, represent spatial sensing demands through importance functions, commonly employ Voronoi partitions \cite{Cortes2004Coverage}, and can be extended to moving regions of interest through time-varying or diffusing importance functions~\cite{Lee2015TimeVaryingDensity,Miah2017NonAutonomousCoverage}.

Despite their effectiveness, moving targets require repeated updates to robot--target assignments or position-dependent partitions, which may increase the computational and communication burden as the swarm size grows. For large-scale MASs, mean-field methods offer a scalable alternative for addressing these limitations by describing the swarm at the population level~\cite{Elamvazhuthi2020Survey}. We therefore adopt a mean-field perspective to describe both the swarm and the tracking demand generated by moving targets at the macroscopic level, with the aim of allowing the swarm to redistribute itself autonomously among the targets.

Mean-field game theory was independently developed by Lasry and Lions~\cite{Lasry2007MeanFieldGames} and by Huang et al.~\cite{Huang2006LargePopulationGames}. Many studies of mean-field swarm control focus on matching prescribed, time-invariant spatial distributions~\cite{Li2017SwarmDensity,Eren2017DensityControl,Krishnan2018SelfOrganization}, with theoretical results on controllability~\cite{Elamvazhuthi2019Controllability} and robustness to density-estimation errors and perturbations~\cite{Zheng2022MeanFieldFeedback,Sinigaglia2025RobustDensity}. Time-varying reference densities have also been considered in field estimation~\cite{Zheng2021FieldEstimation}, causal swarm tracking in Wasserstein space~\cite{Emerick2024CausalTracking}, feedforward density-driven control~\cite{Martinez2026FeedforwardDensity}, and multi-target tracking using Gaussian-mixture references and Voronoi-based coordination~\cite{Abdulghafoor2022MultiTargetTracking,Abdulghafoor2023MultiTargetTracking,Abdulghafoor2024MultiTargetTracking}.

In this paper, we formulate multi-target tracking as a time-varying density-tracking problem, derive a closed-form mean-field feedback law, and analyze the resulting density evolution using Kullback–Leibler (KL) divergence.
Specifically, inspired by density-based formulations of multi-target tracking~\cite{Abdulghafoor2022MultiTargetTracking,Abdulghafoor2023MultiTargetTracking,Abdulghafoor2024MultiTargetTracking}, we represent the tracking and sensing demands associated with moving targets by a time-varying \emph{task-demand density} \(\rho_d\), which specifies the desired spatial allocation of the swarm without explicit agent--target assignments. The control objective is to make the swarm's mean-field density \(\rho\) track \(\rho_d\), with the tracking error measured by the KL divergence \(\operatorname{KL}(\rho\|\rho_d)\). Motivated by the Onsager variational principle~\cite{Doi2011Onsager,ChenLiuXu2025Onsager}, we formulate an instantaneous variational problem to derive a feedback law that accounts for the evolution of \(\rho_d\). The main contributions are as follows:

\begin{enumerate}[label=(\roman*),leftmargin=*,topsep=2pt,itemsep=0pt]
    \item We derive a feedback law for tracking a time-varying task-demand density using an Onsager-type variational formulation. For the resulting closed-loop Fokker–Planck equation, we establish well-posedness and exponential KL convergence under ideal feedback.

    \item We analyze feedback based on estimated swarm densities and filtered target velocities. Treating the estimated drift as given, we establish well-posedness of the corresponding density equation and derive a conditional KL error bound in terms of the feedback mismatch, thereby establishing conditional robustness of KL convergence. The same bound also accommodates velocity saturation through an additional error term.
\end{enumerate}

For density-based target tracking, compared with Voronoi-based tracking methods~\cite{Abdulghafoor2022MultiTargetTracking,Abdulghafoor2023MultiTargetTracking,Abdulghafoor2024MultiTargetTracking}, we directly derive a mean-field feedback for the swarm density to track the desired spatial allocation, without constructing Voronoi partitions. Furthermore, unlike time-varying density-tracking approaches based on prescribed error-diffusion dynamics or predictive optimal control~\cite{Zheng2021FieldEstimation, Emerick2024CausalTracking,Martinez2026FeedforwardDensity}, our approach provides an instantaneous Onsager-type variational alternative that yields a closed-form feedback law, for which we establish KL convergence and conditional robustness for tracking moving task-demand densities. Numerical comparisons with the density-tracking controller in~\cite{Zheng2021FieldEstimation} also show that our controller achieves smaller KDE-based tracking errors at comparable applied control energy in the tested moving-target task.

The remainder of the paper is organized as follows. Section~\ref{section2} formulates the tracking problem. Section~\ref{section3} develops the ideal feedback law, establishes well-posedness, and analyzes convergence. Section~\ref{section4} introduces the estimated feedback law and establishes its conditional robustness properties. Section~\ref{sec:simulation} presents the numerical simulation results, and Section~\ref{section6} concludes the paper.

\section{Problem Formulation}\label{section2}

Let \(\Omega\subset\mathbb R^2\) be a bounded convex Lipschitz task domain and set  \(C:=\overline\Omega\). Consider a swarm of \(N\) homogeneous agents, with the position of agent \(i\) denoted by \(x_i(t)\in C\). In the mean-field description, the swarm is represented by a probability density \(\rho(x,t)\) satisfying
\begin{equation}
    \rho(x,t)\geq0,
    \qquad
    \int_\Omega\rho(x,t)\,\mathrm dx=1.
\end{equation}

To keep the agents within \(C\), we model their microscopic dynamics by a normally reflected stochastic differential equation~\cite{Pilipenko2014}. Since \(\partial \Omega\) is not necessarily smooth, we describe reflection using the outward normal cone
\begin{equation}
    N_C(x):=\left\{q\in\mathbb R^2: \langle q,y-x\rangle\leq0\ \text{for all }y\in C\right\}. \label{eq:agent_normal_cone}
\end{equation}
With diffusion coefficient \(D>0\) and initial position \(x_i(0)\in C\), the dynamics of agent \(i\) are
\begin{equation}
    \mathrm dx_i(t) =u(x_i(t),t)\,\mathrm dt +\sqrt{2D}\,\mathrm dW_i(t) -\nu_i(t)\,\mathrm dL_i(t), \label{eq:reflected_agent_sde}
\end{equation}
where \(u\) is the control drift field and \(W_i\) are independent two-dimensional standard Wiener processes. 

\begin{remark}
    Here, the processes \(x_i\) and \(L_i\) are continuous and adapted, \(L_i\) is nondecreasing with \(L_i(0)=0\). They satisfy, almost surely,
\begin{equation}
    x_i(t)\in C,\qquad
    \int_0^t\mathbf 1_\Omega(x_i(s))\,\mathrm dL_i(s)=0,
    \qquad t\geq0.
    \label{eq:agent_reflection_support}
\end{equation}
The reflection direction \(\nu_i\) is progressively measurable and satisfies
\begin{equation}
    \nu_i(t)\in N_C(x_i(t)),\qquad |\nu_i(t)|=1,\qquad \mathrm dL_i\text{-a.e.} \label{eq:agent_reflection_direction}
\end{equation}
Thus, \(L_i\) increases only when the agent is on the boundary, and the term \(-\nu_i\,\mathrm dL_i\) acts inward to keep its position in \(C\). At a smooth boundary point, we have \(\nu_i = n(x_i)\), where \(n\) is the outward unit normal. 
\end{remark}

For a sufficiently regular density \(\rho\) and \(u\), the density satisfies the following Fokker--Planck (FP) system:
\begin{subequations}\label{eq:rho_u_FP}
\begin{align}
    \partial_t\rho=-\nabla\cdot(\rho u)+D\Delta\rho,\quad &\text{in } \Omega \times (0,T)\label{eq:fp_main},\\
    n\cdot J=n\cdot(\rho u-D\nabla\rho)=0,\quad  &\text{on }\partial\Omega \times (0,T)\label{eq:swarm_no_flux},\\
    \rho(\cdot,0) = \rho_0,\quad & \text{in } \Omega. \label{eq:swarm_initial}
\end{align}
\end{subequations}
Here \(J:=\rho u-D\nabla\rho\) is the total probability current, so \eqref{eq:fp_main} is equivalent to \(\partial_t\rho+\nabla\cdot J=0\). The no-flux condition \eqref{eq:swarm_no_flux} is induced by normal reflection and understood surface-almost everywhere on \(\partial\Omega\), where \(n\) denotes the outward unit normal. Moreover, whenever $\rho>0$, define the current velocity field
\begin{equation}\label{eq:current-velocity}
  v:=\frac{J}{\rho}
  =u-D\nabla\ln\rho.
\end{equation}
Then \eqref{eq:fp_main}--\eqref{eq:swarm_initial}  can be written as 
\begin{equation}\label{eq:rho-current-continuity}
 \begin{cases}
 \partial_t\rho+\nabla\cdot(\rho v)=0,
     &\text{in }\Omega\times(0,T),\\
 n\cdot(\rho v)=0,
     &\text{on }\partial\Omega\times(0,T),\\
 \rho(\cdot,0)=\rho_0,
     &\text{in }\Omega.
 \end{cases}
\end{equation}

The multi-target tracking task is encoded by a desired task-demand density \(\rho_d(x,t)\). Specifically, consider \(M\) targets with positions \(y_m(t)\), \(m=1,\ldots,M\). For each target, there exists a certain demand for surveillance, sensing, or coverage that decreases with the distance from the target. Accordingly, we model this local target-induced demand by a truncated and normalized Gaussian density function centered at \(y_m(t)\):
\begin{align}
    \phi_m(x,t) &= \begin{cases}
        \frac{ \exp\left(-\frac{|x-y_m(t)|^2}{2\sigma_m^2}\right) }{ Z_m(t) }, & x\in \overline{\Omega},\\
        0, & x \notin \overline{\Omega}
    \end{cases}\label{eq:truncated_target_component} \\ 
    Z_m(t) &:= \int_{\Omega} \exp\left(-\frac{|z-y_m(t)|^2}{2\sigma_m^2}\right) \mathrm{d}z.
    \label{eq:truncated_target_gaussian}
\end{align}
Here, \(\sigma_m>0\) specifies the spatial scale of the task demand associated with target \(m\), with larger values indicating wider regions of interest around the target. 

For multiple targets, we aggregate the individual target-induced demands and define the desired \emph{task-demand density}, hereafter referred to as the  \emph{desired density}, by
\begin{equation}
    \rho_d(x,t) := \sum_{m=1}^{M} w_m \phi_m(x,t),\label{eq:task_demand_density}
\end{equation}
where \(w_m\ge 0\) denotes the relative importance of target \(m\), with \(\sum_{m=1}^{M} w_m =1\). These weights specify the prescribed fraction of swarm resources allocated to each target. In the subsequent analysis, $w_m$ and $\sigma_m$ are fixed. 

It's important to note that \(\rho_d\) is not interpreted as a probability density of the targets. Rather, it specifies the desired spatial allocation of the swarm. In particular, for any measurable \(A\subseteq \Omega\), \(\int_A\rho_d(x,t)\,\mathrm{d}x\) represents the desired fraction of agents assigned to \(A\) at time \(t\). Therefore, when \(\rho = \rho_d\), the swarm is considered to have achieved the cooperative multi-target tracking objective.

Accordingly, the tracking performance is measured by the discrepancy between \(\rho\) and \(\rho_d\). We quantify this using the Kullback--Leibler (KL) divergence
\begin{equation}
    F(t) := \mathrm{KL}\bigl(\rho(\cdot,t)\|\rho_d(\cdot,t)\bigr) = \int_\Omega \rho(x,t) \ln\frac{\rho(x,t)}{\rho_d(x,t)} \,\mathrm{d}x. \label{eq:kl_main}
\end{equation}
The control objective is to design \(u(x,t)\) so that the solution \(\rho(x,t)\) of \eqref{eq:fp_main} tracks \(\rho_d(x,t)\) in the KL divergence sense.

\section[Ideal Feedback Law]{Ideal Feedback Law: Design and Convergence}\label{section3}

Throughout this section, unless otherwise stated, \(\rho_d(\cdot,t)\) denotes the prescribed task-demand density in \eqref{eq:task_demand_density}, with prescribed target trajectories satisfying \(y_m\in C^1([0,\infty);\mathbb R^2)\) and \(y_m(t)\in\overline\Omega\) for all \(t\ge0\). Thus, \(\rho_d(\cdot,t)\) is a probability density on \(\Omega\) with a strictly positive extension to \(\overline\Omega\), satisfying
\[
    \rho_d\in C^1(\overline\Omega\times[0,\infty)).
\]
To account for its evolution, we introduce a prescribed velocity field \(v_d\in C([0,\infty);C^1(\overline\Omega;\mathbb R^2))\) satisfying \eqref{eq:desired_transport_and_no_flux}. For static target positions, we take \(v_d\equiv0\).

This section develops the feedback laws in three steps. We first derive conditional KL identities under explicit classical regularity and transport assumptions. This identity motivates an instantaneous Onsager variational problem, whose minimization at fixed densities gives candidate control drift fields. We then establish weak well-posedness of the resulting closed-loop equation for the swarm density \(\rho\), and derive KL dissipation and exponential convergence for its weak solution.

\subsection{Feedback Control Law Design}\label{subsec:feedback_control_law_design}

\paragraph{Conditional KL identity} Assume that the swarm density \(\rho(\cdot,t)\) satisfies \(\rho\in C^1(\overline\Omega\times[0,\infty))\) and \(\rho>0\) on \(\overline\Omega\times[0,\infty)\). Suppose there exists \(v\in C([0,\infty);C^1(\overline\Omega;\mathbb R^2))\) such that \((\rho,v)\) satisfies the continuity equation and no-flux condition in \eqref{eq:rho-current-continuity}, pointwise in \(\Omega\) and surface-almost everywhere on \(\partial\Omega\), respectively, for \(t>0\). The identities below are conditional on these assumptions. 

To account for the evolution of \(\rho_d = \rho_d(\cdot,t)\), suppose that there exists a prescribed transport velocity \(v_d\in C([0,\infty);C^1(\overline\Omega;\mathbb R^2))\) satisfying
\begin{subequations}\label{eq:desired_transport_and_no_flux}
\begin{align}
    \partial_t\rho_d+\nabla\cdot(\rho_dv_d)&=0,
    &&\text{in }\Omega,\quad t>0,
    \label{eq:desired_transport}\\
    n\cdot(\rho_dv_d)&=0,
    &&\text{on }\partial\Omega,\quad t>0,
    \label{eq:desired_no_flux}
\end{align}
\end{subequations}
where the boundary condition is understood surface-almost everywhere. The existence of such a velocity is assumed separately from \eqref{eq:task_demand_density}, and we fix one admissible \(v_d\) without requiring uniqueness. Since \(\rho_d>0\) on \(\overline\Omega\), \eqref{eq:desired_no_flux} implies
\begin{equation}
    n\cdot v_d=0
    \quad\text{surface-a.e. on }\partial\Omega,\quad t>0.
    \label{eq:vd_tangent}
\end{equation}
For a static desired density, we take \(v_d=0\). In the nominal analysis, \(v_d\) always denotes this reference transport velocity. 

Differentiating \eqref{eq:kl_main} gives
\begin{align}
    \frac{\mathrm{d}}{\mathrm{d}t}F(t)
    = \int_{\Omega} \left(1+\ln\frac{\rho}{\rho_d}\right) \partial_t\rho\,\mathrm{d}x - \int_{\Omega} \frac{\rho}{\rho_d}\partial_t\rho_d\,\mathrm{d}x.
    \label{eq:entropy-derivative-first}
\end{align}
Substituting the swarm continuity equation and \eqref{eq:desired_transport} into \eqref{eq:entropy-derivative-first}, then integrating by parts using the no-flux conditions, gives 
\begin{equation}
     \int_{\Omega} \left(1+\ln\frac{\rho}{\rho_d}\right) \partial_t\rho\,\mathrm{d}x = \int_{\Omega} \rho v\cdot\nabla\ln\frac{\rho}{\rho_d}\,\mathrm{d}x. \label{eq:first-transport-term}
\end{equation}
\begin{align}
    \int_{\Omega}
    \frac{\rho}{\rho_d}\,
    \nabla\cdot(\rho_dv_d)\,\mathrm{d}x
    ={}&
    -\int_{\Omega}
    \rho_dv_d\cdot
    \nabla\left(\frac{\rho}{\rho_d}\right)\,\mathrm{d}x
    \notag\\
    ={}&
    -\int_{\Omega}
    \rho v_d\cdot
    \nabla\ln\frac{\rho}{\rho_d}\,\mathrm{d}x,
    \label{eq:second-transport-term}
\end{align}
Combining \eqref{eq:first-transport-term} and
\eqref{eq:second-transport-term}, we obtain, for every \(t>0\),
\begin{equation}
    \frac{\mathrm{d}}{\mathrm{d}t}F(t) = \int_{\Omega} \rho\nabla\ln\frac{\rho}{\rho_d}\cdot(v-v_d) \,\mathrm{d}x. \label{eq:entropy-identity}
\end{equation}
Equation \eqref{eq:entropy-identity} shows that the instantaneous change of the tracking error is governed by \(\xi:=v-v_d\), which we call the relative current velocity field.

\paragraph{Onsager Selection of the Feedback Laws}

Fix \(t\ge0\) and a probability density \(\rho(\cdot,t)\in C^1(\overline\Omega)\) (no need to satisfy any evolution equation), strictly positive on \(\overline\Omega\). During the minimization, \(\rho(\cdot,t)\), \(\rho_d(\cdot,t)\), and \(v_d(\cdot,t)\) are held fixed. For vector fields \(z\in L^2(\Omega,\rho\,\mathrm{d}x;\mathbb R^2)\), define
\begin{equation}
    \|z\|_\rho^2:=\int_\Omega\rho\|z\|^2\,\mathrm{d}x.
    \label{eq:weighted_velocity_norm}
\end{equation}

Inspired by \eqref{eq:entropy-identity}, where \(\dot{F}\) is governed by the relative current velocity, define
\[
    G(\xi) := \int_\Omega \rho \nabla\ln\frac{\rho}{\rho_d} \cdot \xi\,\mathrm{d}x,
\]
where \(\xi\in L^2(\Omega,\rho\mathrm{d}x;\mathbb R^2)\) represents a trial relative current velocity. When $\nabla\ln\frac{\rho}{\rho_d}\ne0$, minimizing \(G(\xi)\) alone is unbounded below, since the magnitude of a descent direction can be increased arbitrarily. Motivated by the Onsager principle \cite{ChenLiuXu2025Onsager,Doi2011Onsager}, we combine $G(\xi)$ with a quadratic transport cost to define the Rayleighian
\begin{equation}
    \Psi_\rho(\xi) := \frac{1}{2K} \int_\Omega\rho\|\xi\|^2\,\mathrm{d}x + \int_\Omega \rho \nabla\ln\frac{\rho}{\rho_d}\cdot \xi\,\mathrm{d}x, \label{eq:rayleighian_1}
\end{equation}
where \(K > 0\) is a prescribed gain coefficient. We then select a candidate current velocity by solving the Onsager variational problem
\[
    \min_{\xi\in L^2(\Omega,\rho\,\mathrm{d}x;\mathbb R^2)} \Psi_\rho(\xi).
\]
Note that \(\Psi_\rho(\xi) = \frac{1}{2K} \left\| \xi+K \nabla\ln\frac{\rho}{\rho_d} \right\|_\rho^2 - \frac{K}{2} \left\| \nabla\ln\frac{\rho}{\rho_d} \right\|_\rho^2\), thus the unique minimizer \(-K \nabla\ln\frac{\rho}{\rho_d}\) gives the candidate relative current velocity \(\xi^* = -K\nabla\ln\frac{\rho}{\rho_d}\).  And by \(\xi = v - v_d\) and \(u = v + D\nabla\ln\rho \), the candidate control law is
\begin{equation}
    \boxed{u^\star = v_d+K\nabla\ln\rho_d-(K-D)\nabla\ln\rho.}\label{eq:diffusive_time_varying_control}
\end{equation}

We adopt \eqref{eq:diffusive_time_varying_control} as the feedback law for the prescribed reference pair \((\rho_d,v_d)\). Formal substitution into \eqref{eq:rho_u_FP}, using \(\rho u^\star-D\nabla\rho=\rho(v_d+K\nabla\ln\rho_d)-K\nabla\rho\), gives, on any finite interval \([0,T]\), the closed-loop system
\begin{equation}
    \begin{cases}
    \partial_t\rho = K\Delta\rho-\nabla\cdot\!\left[\rho\bigl(v_d+K\nabla\ln\rho_d\bigr)\right], & \text{in }\Omega\times(0,T),\\
    n\cdot\!\left[\rho\bigl(v_d+K\nabla\ln\rho_d\bigr)-K\nabla\rho\right]=0, & \text{on }\partial\Omega\times(0,T),\\
    \rho(\cdot,0) = \rho_0(\cdot), & \text{in } \Omega.
\end{cases}\label{eq:closed_loop_system}
\end{equation}
Here \(\rho_d\) is defined in \eqref{eq:task_demand_density} and \(v_d\in C([0,\infty);C^1(\overline\Omega;\mathbb R^2))\) is its prescribed transport velocity satisfying \eqref{eq:desired_transport_and_no_flux}. 

\subsection{Closed-Loop Well-Posedness and KL Dissipation}
\label{subsec:closed_loop_wellposedness}

With \(b:=v_d+K\nabla\ln\rho_d\), the closed-loop system becomes 
\begin{equation}
    \begin{cases}
    \partial_t\rho=K\Delta\rho-\nabla\cdot(\rho b),
        &\text{in }\Omega\times(0,T),\\
    n\cdot(\rho b-K\nabla\rho)=0,
        &\text{on }\partial\Omega\times(0,T),\\
    \rho(\cdot,0)=\rho_0,
        &\text{in }\Omega.
    \end{cases}
    \label{eq:closed_loop_wp_problem}
\end{equation}
Multiplying the PDE in \eqref{eq:closed_loop_wp_problem} by $\psi\in V:=H^1(\Omega)$ and formally integrating by parts, the no-flux boundary condition gives
\[
    \langle\partial_t\rho,\psi\rangle
    =-K\int_\Omega\nabla\rho\cdot\nabla\psi\,\mathrm{d}x
      +\int_\Omega\rho b\cdot\nabla\psi\,\mathrm{d}x.
\]
Accordingly, a function $\rho\in L^2(0,T;V)\cap H^1(0,T;V^*)$ is a weak solution of \eqref{eq:closed_loop_wp_problem} if it satisfies this identity for every $\psi\in V$ and a.e. $t\in(0,T)$, with $\rho(0)=\rho_0$ in $L^2(\Omega)$. 

\begin{proposition}[Closed-loop well-posedness]
\label{prop:closed_loop_wellposedness}
Let \(\Omega\) be the task domain of Section~\ref{section2}, and fix \(K>0\), \(T>0\), and \(D\ge0\). Under the assumptions of this section, the prescribed reference pair \((\rho_d,v_d)\) satisfies
\begin{subequations}\label{eq:closed_loop_wp_data}
\begin{align}
    &\rho_d\in L^\infty(0,T;W^{1,\infty}(\Omega)),
    \label{eq:closed_loop_wp_data_L_infty}\\
    &0<\ell_d\le\rho_d(\cdot,t)\le L_d<\infty
    \quad\text{a.e. in }\Omega,\quad t\ge0,
    \label{eq:closed_loop_wp_data_boundedness}\\
    &v_d\in L^\infty(\Omega\times(0,T);\mathbb R^2),
    \label{eq:closed_loop_vd_data_L_infty}\\
    &\rho_d\in C^1([0,T];W^{1,\infty}(\Omega)).
    \label{eq:reference_time_regularity}
\end{align}
\end{subequations}
If \(\rho_0\in L^2(\Omega)\), \(\rho_0\ge0\) a.e. in \(\Omega\), and \(\int_\Omega\rho_0\,\mathrm{d}x=1\), then the following properties hold: 
\begin{enumerate}
    \item \eqref{eq:closed_loop_wp_problem} admits a unique
    weak solution
\begin{equation}
    \rho\in C([0,T];L^2(\Omega))
    \cap L^2(0,T;H^1(\Omega)),
    \label{eq:closed_loop_wp_spaces}
\end{equation}
with \(\partial_t\rho\in L^2(0,T;(H^1(\Omega))^*)\), \(\rho(0)=\rho_0\) in \(L^2(\Omega)\), and
\begin{equation}
    \begin{aligned}
    \langle\partial_t\rho,\psi\rangle
    &=-K\int_\Omega\nabla\rho\cdot\nabla\psi\,\mathrm{d}x\\
    &\quad+\int_\Omega\rho b\cdot\nabla\psi\,\mathrm{d}x,
    \qquad\psi\in H^1(\Omega),
    \end{aligned}
    \label{eq:closed_loop_wp_weak}
\end{equation}
for a.e. \(t\in(0,T)\). Here \(\langle\cdot,\cdot\rangle\) denotes the duality pairing between \((H^1(\Omega))^*\) and \(H^1(\Omega)\).

    \item For every \(t\in[0,T]\), the solution satisfies \(\rho(\cdot,t)\ge0\) a.e. in \(\Omega\) and \(\int_\Omega\rho(x,t)\,\mathrm{d}x=1\). For two initial densities \(\rho_{0,1}\) and \(\rho_{0,2}\)), their solutions \(\rho_1\) and \(\rho_2\) satisfy
    
\begin{equation}
\begin{gathered}
    \|\rho_1(t)-\rho_2(t)\|_{L^2(\Omega)}^2 \le e^{\frac{B_T^2t}{K}} \|\rho_{0,1}-\rho_{0,2}\|_{L^2(\Omega)}^2, \\
    B_T:=\|b\|_{L^\infty(\Omega\times(0,T);\mathbb R^2)}.
\end{gathered}
    \label{eq:closed_loop_continuous_dependence}
\end{equation}

    \item The ratio \(q:=\rho/\rho_d\) satisfies
\begin{equation}
    \begin{gathered}
    q\in C([0,T];L^2(\Omega))
    \cap L^2(0,T;H^1(\Omega)),\\
    \partial_tq\in L^2(0,T;(H^1(\Omega))^*).
    \end{gathered}
    \label{eq:closed_loop_ratio_spaces}
\end{equation}

\item If, additionally, \(\rho_0>0\) a.e. in \(\Omega\), then, for every \(t\in[0,T]\), \(\rho(\cdot,t)>0\) a.e. in \(\Omega\).
\end{enumerate}
Moreover, since \(T>0\) is arbitrary , uniqueness yields a global weak solution.
\end{proposition}

\begin{proof}
1. Set $Q_T:=\Omega\times(0,T)$, $V:=H^1(\Omega)$, and $H:=L^2(\Omega)$. The bounds in
\eqref{eq:closed_loop_wp_data} give \(B_T:=\|b\|_{L^\infty(Q_T)} \leq\|v_d\|_{L^\infty(Q_T)}  +\frac{K}{\ell_d}\|\nabla\rho_d\|_{L^\infty(Q_T)}<\infty\). Define
\[
    a_t(f,g):=K\int_\Omega\nabla f\cdot\nabla g\,\mathrm{d}x
       -\int_\Omega f b\cdot\nabla g\,\mathrm{d}x,
    \qquad f,g\in V.
\]
For each $f,g\in V$, the map $t\mapsto a_t(f,g)$ is measurable. Writing $\|\cdot\|_2$ for the $L^2(\Omega)$ norm, the Cauchy--Schwarz inequality gives
\[
    \begin{aligned}
    |a_t(f,g)| & \leq K\|\nabla f\|_2\|\nabla g\|_2 +B_T\|f\|_2\|\nabla g\|_2\\
    &\leq(K+B_T)\|f\|_V\|g\|_V.
    \end{aligned}
\]
By Young's inequality, \(a_t(f,f)\geq\frac K2\|\nabla f\|_2^2-\frac{B_T^2}{2K}\|f\|_2^2\). Consequently, with $\omega_T:=B_T^2/(2K)+K/2$,
\[
    a_t(f,f)+\omega_T\|f\|_H^2\geq\frac K2\|f\|_V^2.
\]
Hence $a$ is a non-autonomous closed form on \(H\) \cite{ArendtDierOuhabaz2014}. It defines $A(t)\in\mathcal L(V,V^*)$ through $\langle A(t)f,g\rangle=a_t(f,g)$. Since $V$ is densely and continuously embedded in $H$, Lions' theorem \cite[Section 1]{ArendtDierOuhabaz2014} yields a unique solution $\rho\in L^2(0,T;V)\cap H^1(0,T;V^*)$ of
\[
    \partial_t\rho+A(t)\rho=0\quad\text{in }V^*,\quad \text{a.e. } t\in [0,T];
    \qquad\rho(0)=\rho_0.
\]
By the definition of \(A(t)\), this equation is equivalent to \eqref{eq:closed_loop_wp_weak}. It is well known that \(L^2(0,T;V)\cap H^1(0,T;V^*) \hookrightarrow C([0,T];H)\) and hence \(\rho\) has a unique continuous representative, also denoted by \(\rho\). Thus the initial condition makes sense in \(H = L^2(\Omega)\), and this \(\rho \in C([0,T];H)\cap L^2(0,T;V)\) is a unique weak solution of \eqref{eq:closed_loop_wp_problem}. Moreover, \(\rho \in H^1(0,T;V^*)\) implies \(\partial_t\rho\in L^2(0,T;V^*)\).

2. Set \(\rho^-:=\max\{-\rho,0\}\). Testing \eqref{eq:closed_loop_wp_weak} with \(-\rho^-\) and using the truncation chain rule gives, for a.e.
\(t\in(0,T)\),
\[
    \begin{aligned}
    \frac12\frac{\mathrm{d}}{\mathrm{d}t}\|\rho^-\|_2^2
        +K\|\nabla\rho^-\|_2^2
    &=\int_\Omega\rho^-b\cdot\nabla\rho^-\,\mathrm{d}x\\
    &\leq\frac K2\|\nabla\rho^-\|_2^2 + \frac{B_T^2}{2K}\|\rho^-\|_2^2.
    \end{aligned}
\]
Gronwall's inequality therefore yields \(\|\rho^-\|_2^2 \leq e^{B_T^2t/K}\|\rho^-_0\|_2^2 = 0\), so $\rho^-=0$ in \(L^2(\Omega)\) for every \(t\in [0,T]\), and thus \(\rho(\cdot,t)\ge 0\) a.e. in \(\Omega\) for every \(t\in [0,T]\). Next, taking \(\psi\equiv1\) in \eqref{eq:closed_loop_wp_weak} shows that the absolutely continuous function \(t\mapsto\int_\Omega\rho(x,t)\,\mathrm{d}x\) has zero derivative. Together with the initial condition and the continuous \(L^2\) representative, this gives \(\int_\Omega\rho(x,t)\,\mathrm{d}x=1\) for every \(t\in[0,T]\). The same energy calculation for \(\rho_1-\rho_2\) proves \eqref{eq:closed_loop_continuous_dependence}
for every \(t\in[0,T]\).

3. Set
$a:=1/\rho_d$. The regularity
\eqref{eq:reference_time_regularity} and the lower bound
\(\rho_d\geq\ell_d>0\) imply $a\in C^1([0,T];W^{1,\infty}(\Omega))$.
By the spatial product rule, there exists a constant $C_T>0$, independent of $t$ and $f$, such that \(\|a(t)f\|_V\leq C_T\|f\|_V, \ \forall t\in[0,T],\ f\in V\). Thus, for $q:=a\rho$,
\[
    \int_0^T\|q(t)\|_V^2\,\mathrm{d}t
    \leq C_T^2\int_0^T\|\rho(t)\|_V^2\,\mathrm{d}t<\infty.
\]
Moreover, writing \(q(t)-q(s) =a(t)\bigl(\rho(t)-\rho(s)\bigr)  +\bigl(a(t)-a(s)\bigr)\rho(s)\), we obtain
\[
\begin{aligned}
    \|q(t)-q(s)\|_H
    &\leq \ell_d^{-1}\|\rho(t)-\rho(s)\|_H\\
    &\quad+\|a(t)-a(s)\|_{L^\infty(\Omega)}\|\rho(s)\|_H.
\end{aligned}
\]
Since $\rho\in C([0,T];H)$ and $a\in C([0,T];L^\infty(\Omega))$, the right-hand side tends to \(0\) as $t\to s$. Hence \(q=\frac{\rho}{\rho_d} \in C([0,T];H)\cap L^2(0,T;V)\).

For every $\psi\in V$, the weak product rule gives 
\[
    \langle\partial_tq,\psi\rangle =\left\langle\partial_t\rho,\frac{\psi}{\rho_d}\right\rangle - \int_\Omega\frac{\rho\,\partial_t\rho_d}{\rho_d^2}\psi\,\mathrm{d}x.
\]
Since  $\partial_t\rho_d/\rho_d^2\in L^\infty(Q_T)$, set \(\Tilde{C}_T = \max \left\{C_T,\left\|{\partial_t\rho_d}/{\rho_d^2}\right\|_{L^\infty(Q_T)}\right\} \), we obtain \(\left|\left\langle\partial_t\rho,\frac{\psi}{\rho_d}\right\rangle\right|\leq \Tilde{C}_T\|\partial_t\rho\|_{V^*}\|\psi\|_V\) and \(\left|\int_\Omega \frac{\rho\,\partial_t\rho_d}{\rho_d^2} \psi\,\mathrm{d}x\right|\leq \Tilde{C}_T\|\rho\|_2\|\psi\|_V\). Therefore, 
\[
    \|\partial_tq\|_{V^*}  =\sup_{\|\psi\|_V\leq1}|\langle\partial_tq,\psi\rangle|\leq \Tilde{C}_T\bigl(\|\partial_t\rho\|_{V^*}+\|\rho\|_2\bigr).
\]
Since $\partial_t\rho\in L^2(0,T;V^*)$ and
$\rho\in L^2(0,T;H)$, it follows that
$\partial_tq\in L^2(0,T;V^*)$. Together with the regularity of $q$ established above, this proves \eqref{eq:closed_loop_ratio_spaces}.

4. Suppose \(\rho_0>0\) a.e. in \(\Omega\). Since \(b=v_d+K\nabla\ln\rho_d\), \(\rho_d\) satisfies the same linear weak equation as \(\rho\):
\[
    \langle\partial_t f,\psi\rangle =-K\int_\Omega\nabla f\cdot\nabla\psi\,dx   +\int_\Omega fb\cdot\nabla\psi\,dx, \qquad \psi\in H^1(\Omega),
\]
for \(f=\rho_d\) and a.e. \(t\in(0,T)\), with initial value \(\rho_d(\cdot,0)\). For each integer \(n\ge1\), define \(f_{n,0}:=\min\{\rho_d(\cdot,0),n\rho_0\}\), and let \(f_n\) solve the same linear no-flux problem with initial value \(f_{n,0}\). Since preceding arguments do not require initial data to have unit mass, they also apply to \(f_{n,0}\).

Because \(\rho_0>0\) a.e., we have \(f_{n,0}\to \rho_d(\cdot,0)\), a.e. in \(\Omega\). Moreover, \(0\le f_{n,0}\le\rho_d(\cdot,0)\), so dominated convergence gives \(\|f_{n,0}-\rho_d(\cdot,0)\|_2\to 0\). Applying the continuous-dependence estimate then yields
\[
    \sup_{t\in[0,T]}\|f_n(t)-\rho_d(t)\|_2^2
    \le e^{B_T^2T/K}
       \|f_{n,0}-\rho_d(\cdot,0)\|_2^2
    \to 0.
\]
Since \(n\rho-f_n\) satisfies the same equation and has nonnegative initial value, the negative-part argument used in Step 2 gives
\[
    0\le f_n(t)\le n\rho(t)
    \quad\text{a.e. in }\Omega,
    \qquad t\in[0,T].
\]
Fix any \(t\in[0,T]\) and set
\(A_t:=\{x\in\Omega:\rho(x,t)=0\}\).
Then \(f_n(t)=0\) a.e. on \(A_t\), and consequently
\[
    \ell_d^2|A_t|
    \le\int_{A_t}\rho_d(x,t)^2\,\mathrm{d}x
    \le\|\rho_d(t)-f_n(t)\|_2^2
    \longrightarrow0.
\]
Thus \(|A_t|=0\). Since \(t\) was arbitrary,
\(\rho(\cdot,t)>0\) a.e. in \(\Omega\) for every \(t\in[0,T]\).
\end{proof}

\begin{remark}
    The regularity in \eqref{eq:closed_loop_wp_data} follows from the standing assumptions. The uniform bounds on \(\rho_d\) are obtained as follows. Let \(d_\Omega:=\operatorname{diam}(\Omega)\) and \(a_m:=d_\Omega^2/(2\sigma_m^2)\). Since \(x,y_m(t)\in\overline\Omega\),
\[
    |\Omega|e^{-a_m}\le Z_m(t)\le|\Omega|,
    \qquad
    \frac{e^{-a_m}}{|\Omega|}\le\phi_m(x,t)
    \le\frac{e^{a_m}}{|\Omega|}.
\]
Consequently, with \(\ell_d:=\min_m\frac{e^{-a_m}}{|\Omega|}\) and \(L_d:=\max_m\frac{e^{a_m}}{|\Omega|}\), 
\begin{equation}
    0<\ell_d\le\rho_d(x,t)\le L_d<\infty \label{eq:uniform_desired_density_bounds}
\end{equation}
for all \((x,t)\in\overline\Omega\times[0,\infty)\). Thus \eqref{eq:closed_loop_wp_data_boundedness} holds.
\end{remark}

Define the relative Fisher information:
\begin{equation}
    \mathcal I(\rho\|\rho_d)
    :=4\int_\Omega\rho_d
       \left\|\nabla\sqrt{\frac{\rho}{\rho_d}}\right\|^2
       \,\mathrm{d}x,
    \label{eq:relative_fisher_information}
\end{equation}
when $\sqrt{\rho/\rho_d}\in H^1(\Omega)$, and set
$\mathcal I(\rho\|\rho_d)=+\infty$ otherwise. For positive,
sufficiently regular densities, \eqref{eq:relative_fisher_information} agrees with $\int_\Omega\rho\|\nabla\ln(\rho/\rho_d)\|^2\,\mathrm{d}x$.

\begin{proposition}[KL dissipation]
\label{prop:KL_divergence_derivative}
Fix \(T>0\), and let \(\rho\) be the weak solution constructed in Proposition~\ref{prop:closed_loop_wellposedness} on \([0,T]\). Then \(\sqrt{\rho/\rho_d}\in L^2(0,T;H^1(\Omega))\) and \(F(t):=\operatorname{KL}(\rho(\cdot,t)\|\rho_d(\cdot,t)) \in W^{1,1}(0,T)\). Moreover, for every \(0\le s\le t\le T\),
\begin{equation}
    F(t)+K\int_s^t \mathcal I(\rho(\tau)\|\rho_d(\tau))\,\mathrm{d}\tau = F(s). \label{eq:integral_kl_dissipation}
\end{equation}
In particular,
\begin{equation}
    F'(t)=-K\mathcal I(\rho(t)\|\rho_d(t))\leq0 \quad\text{for a.e. }t\in(0,T). \label{eq:time_varying_kl_dissipation}
\end{equation}
\end{proposition}

\begin{proof}
Write \(q:=\rho/\rho_d\) and \(b:=v_d+K\nabla\ln\rho_d\). The classical reference transport and no-flux conditions \eqref{eq:desired_transport_and_no_flux} imply
\begin{equation}
    \int_\Omega\partial_t\rho_d\,\psi\,\mathrm{d}x =\int_\Omega\rho_d v_d\cdot\nabla\psi\,\mathrm{d}x, \qquad\psi\in H^1(\Omega), \label{eq:reference_transport_weak}
\end{equation}
for a.e. \(t\in(0,T)\). Since \(\rho_d b-K\nabla\rho_d=\rho_dv_d\), it gives
\begin{equation}
    \langle\partial_t f,\psi\rangle =-K\int_\Omega\nabla f\cdot\nabla\psi\,dx +\int_\Omega fb\cdot\nabla\psi\,dx, \quad \psi\in H^1(\Omega),\label{eq:same_linear_parabolic_equation}
\end{equation}
for a.e. \(t\in(0,T)\), where \(f = \rho_d\). By Proposition \ref{prop:closed_loop_wellposedness}, \(\rho\) satisfies \eqref{eq:same_linear_parabolic_equation} with \(f\) replaced by \(\rho\). Thus \(\rho_d\) and \(\rho\) solve the same linear parabolic equation in the weak sense.

For \(0<\varepsilon\le1\), set
\[
    f_\varepsilon:=\rho+\varepsilon \rho_d, \qquad q_\varepsilon:=\frac{f_\varepsilon}{\rho_d}=q+\varepsilon.
\]
By linearity, \(f_\varepsilon\) satisfies \eqref{eq:same_linear_parabolic_equation} with \(f = f_{\varepsilon}\), and its associated current is \(J_\varepsilon =f_\varepsilon b-K\nabla f_\varepsilon =\rho_d q_\varepsilon v_d-K\rho_d \nabla q\). Now define
\[
    E_\varepsilon(t) :=\int_\Omega \rho_d\, \eta (q_{\varepsilon}) \,dx,
\]
where \(\eta(s)=s\ln s-s+1 \), then by \(\rho_dq_\varepsilon=f_\varepsilon\) we obtain
\[
E_\varepsilon =\int_\Omega\eta(f_\varepsilon)\,dx -\int_\Omega f_\varepsilon\ln\rho_d\,dx +\int_\Omega(\rho_d-1)\,dx.
\]
Since \(f_\varepsilon=\rho+\varepsilon \rho_d\ge\varepsilon \ell_d>0\), the function \(\eta(s)\) can be extended below \(\varepsilon \ell_d\) to a \(C^2\) function on \(\mathbb R\) with bounded second derivative, without changing its values along \(f_\varepsilon\). Applying the variational chain rule \cite[Lemma~A.19]{Jungel2025CrossDiffusionNotes} to \(f_\varepsilon\), together with the weak product rule for \(f_\varepsilon\ln \rho_d\), gives \(E_\varepsilon\in W^{1,1}(0,T)\) and
\[
\begin{aligned}
E_\varepsilon'
&=\langle\partial_t f_\varepsilon,\ln f_\varepsilon\rangle
-\langle\partial_t f_\varepsilon,\ln \rho_d\rangle
-\int_\Omega\frac{f_\varepsilon}{ \rho_d }\,\partial_t \rho_d\,dx
+\int_\Omega\partial_t \rho_d\,dx\\
&=\langle\partial_t f_\varepsilon,\ln q_\varepsilon\rangle
+\int_\Omega(1-q_\varepsilon)\,\partial_t \rho_d\,dx
\end{aligned}
\]
for a.e. \(t\in(0,T)\). The product rule applies because \(\ln \rho_d\in C^1([0,T];W^{1,\infty}(\Omega))\). Testing \eqref{eq:same_linear_parabolic_equation} for \(f=f_\varepsilon\) with \(\ln q_\varepsilon\), and \eqref{eq:reference_transport_weak} with \(1-q_\varepsilon\), gives
\[
\begin{aligned}
E_\varepsilon'
&=\int_\Omega
 \bigl(\rho_d q_\varepsilon v_d-K \rho_d \nabla q\bigr)
 \cdot\frac{\nabla q}{q_\varepsilon}\,dx
 -\int_\Omega  \rho_dv_d\cdot\nabla q\,dx\\
&=-K\int_\Omega  \rho_d\frac{|\nabla q|^2}{q+\varepsilon}\,dx.
\end{aligned}
\]
Thus, for every \(0\le s\le t\le T\),
\begin{equation}
    E_\varepsilon(t) + K\int_s^t\!\int_\Omega \rho_d\frac{|\nabla q|^2}{q+\varepsilon}\,dx\,d\tau = E_\varepsilon(s).\label{eq:integral_of_E}
\end{equation}

Fix \(t\in[0,T]\) and suppress the time argument. Define \(0\ln0:=0\), then \(\eta\) is continuous on \([0,\infty)\), and then
\[
    \rho_d(x)\eta(q(x)+\varepsilon)\longrightarrow \rho_d(x)\eta(q(x)) \quad\text{for a.e. }x\in\Omega
\]
as \(\varepsilon\downarrow 0\). Observe that \(0\le \eta(r)\le(r-1)^2\ (r\ge 0)\), then for every \(0<\varepsilon\le1\),
\[
    0\le \rho_d\,\eta(q+\varepsilon) \le L_d(q+\varepsilon-1)^2\le 2L_d(q^2+1).
\]
The function \(2L_d(q^2+1)\) is independent of \(\varepsilon\) and belongs to \(L^1(\Omega)\), because \(q(t)\in L^2(\Omega)\) and \(|\Omega|<\infty\). The dominated convergence theorem therefore gives
\[
    \lim_{\varepsilon\downarrow0}E_\varepsilon(t) = \int_\Omega \rho_d\,\eta(q)\,dx.
\]
Using \(\rho_d q=\rho\) and \(\int_\Omega\rho\,dx=\int_\Omega \rho_d \,dx=1\), we obtain
\[    
    \int_\Omega \rho_d\,\eta(q)\,dx = \int_\Omega\rho\ln\frac{\rho}{\rho_d}\,dx -\int_\Omega\rho\,dx+\int_\Omega \rho_d\,dx=F(t).
\]
Hence \(E_\varepsilon(t)\to F(t)\) for every \(t\in[0,T]\). Since \(\rho_d\frac{|\nabla q|^2}{q+\varepsilon}\) in  \eqref{eq:integral_of_E} increases as \(\varepsilon\downarrow 0\), monotone convergence theorem  yields
\begin{equation}
    F(t)+K\int_s^t\!\int_\Omega \rho_d\frac{|\nabla q|^2}{q}\,dx\,d\tau=F(s),\label{eq:entropy_balance}
\end{equation}
where the quotient is defined to be zero on \(\{q=0\}\).

Taking \(s=0\) and \(t=T\) in \eqref{eq:entropy_balance} gives \(\int_0^T\!\int_\Omega \rho_d\frac{|\nabla q|^2}{q}\,dx\,dt = \frac{F(0) - F(T)}{K}\). Together with \(\rho_d\ge \ell_d>0\) and \(q\in L^1(\Omega\times(0,T))\), we obtain a uniform bound for \(\sqrt{q+\varepsilon}\) in \(L^2(0,T;H^1(\Omega))\):
\begin{align*}
    & {}\|\sqrt{q+\varepsilon}\|_{L^2(0,T;H^1(\Omega))}^2 = \int_{0}^T \int_{\Omega}\left(|\sqrt{q+\varepsilon}|^2 + |\nabla \sqrt{q+\varepsilon}|^2\right)\mathrm{d}x \mathrm{d}t \\
    & = \int_{0}^T \int_{\Omega} \left(q + \varepsilon + \frac{|\nabla q|^2}{4(q+\varepsilon)}\right)\mathrm{d}x \mathrm{d}t \\
    & \le \|q\|_{L^1(Q_T)} + T \ |\Omega| + \int_{0}^T \int_{\Omega}\frac{|\nabla q|^2}{4q}\cdot \frac{\rho_d}{\ell_d}\mathrm{d}x \mathrm{d}t \\
    & = \|q\|_{L^1(Q_T)} + T \ |\Omega| + \frac{F(0) - F(T)}{4\ell_dK} < \infty,
\end{align*} 
where \(Q_T = \Omega \times (0,T)\).
Since \(\bigl|\sqrt{q+\varepsilon}-\sqrt q\bigr|^2\le\varepsilon\), the family \(\sqrt{q+\varepsilon}\ (\varepsilon>0)\) converges strongly to \(\sqrt q\) in \(L^2(\Omega\times(0,T))\) as \(\varepsilon\downarrow0\). Together with the uniform \(L^2(0,T;H^1(\Omega))\) bound, this implies, by weak compactness in \(L^2(0,T;H^1(\Omega))\), that
\[
\sqrt q\in L^2(0,T;H^1(\Omega)).
\]
The Sobolev chain rule on \(\{q>0\}\), together with the fact that \(\nabla\sqrt q=0\) almost everywhere on \(\{q=0\}\), yields
\[
\frac{|\nabla q|^2}{q}=4|\nabla\sqrt q|^2
\quad\text{for a.e. }(x,t)\in\Omega\times(0,T),
\]
where the quotient is defined to be zero on \(\{q=0\}\). Together with \eqref{eq:entropy_balance} and \eqref{eq:relative_fisher_information}, we obtain
\[
F(t)+K\int_s^t
\mathcal I(\rho(\tau)\|\rho_d(\tau))\,d\tau=F(s).
\]
Taking \(s=0\) and \(t=T\) shows that \(\mathcal I(\rho(\cdot)\|\rho_d(\cdot))\in L^1(0,T)\). And taking \(s = 0\) and \(t \in (0,T]\) gives
\[
F(t)=F(0)-K\int_0^t\mathcal I(\rho(\tau)\|\rho_d(\tau))\,d\tau.
\]
Thus \(F\) is absolutely continuous, belongs to \(W^{1,1}(0,T)\), and satisfies
\(F'=-K\mathcal I(\rho\|\rho_d)\) almost everywhere in \((0,T)\).
\end{proof}

Note that the closed-loop equation can be written as
\[
\partial_t\rho+\nabla\cdot(\rho v_d)
=K\nabla\cdot\left(\rho_d\nabla\frac{\rho}{\rho_d}\right).
\]
When \(\rho=\rho_d\), the ratio \(\rho/\rho_d\) equals one, so the correction term vanishes and the swarm density obeys the same transport equation as the desired density, with velocity field \(v_d\). Thus, \(\rho=\rho_d\) is an exact solution of the closed-loop equation. By uniqueness, if \(\rho(s)=\rho_d(s)\) at some time \(s\), the two densities remain equal and are transported together by \(v_d\) for all subsequent times.

The entropy identity \eqref{eq:integral_kl_dissipation} provides an equivalent verification: \(F(s)=0\), together with the nonnegativity of the KL divergence and its dissipation, implies \(F(t)=0\) for every \(t\ge s\), and hence \(\rho(t)=\rho_d(t)\) almost everywhere in \(\Omega\).

\subsection{Exponential Convergence of the KL Divergence}

\begin{lemma}
\label{lem:uniform_lsi_time_varying}
Let \(\Omega\) be the task domain of Section~\ref{section2}. Under the assumptions of this section, there exists \(C_{\mathrm{KL}}^{\mathrm{unif}}>0\), depending only on \(\Omega,\ell_d,L_d\), such that, for every \(t\ge0\) and
every nonnegative probability density \(\rho\) with \(\sqrt{\rho/\rho_d(\cdot,t)}\in H^1(\Omega)\),
\begin{equation}
    \operatorname{KL}\bigl(\rho\|\rho_d(\cdot,t)\bigr) \le C_{\mathrm{KL}}^{\mathrm{unif}} \mathcal I\bigl(\rho\|\rho_d(\cdot,t)\bigr). \label{eq:uniform_lsi_time_varying}
\end{equation}
\end{lemma}

\begin{proof}
Fix \(t\ge0\), choose \(\kappa>0\), and set \(\psi(x):=\kappa\|x - x_0\|^2/2\) where \(x_0\in \Omega\). Define 
\[
    \mu(dx):=Z_\mu^{-1}e^{-\psi(x)}\,dx,\quad \nu_t(dx):=\rho_d(x,t)\,dx,
\]
where \(Z_\mu:=\int_\Omega e^{-\psi(x)}\,dx\). For a probability measure \(\nu\), write \(\operatorname{Ent}_\nu(f^2) :=\int f^2\ln f^2\,d\nu -\left(\int f^2\,d\nu\right) \ln\left(\int f^2\,d\nu\right)\), with \(0\ln0:=0\). Since \(\Omega\) is convex and \(\nabla^2\psi=\kappa I\), by the Bakry--\'{E}mery criterion~\cite[Theorem~A1]{Schlichting2019Mixtures}, \(\mu\) satisfies:
\[
    \operatorname{Ent}_\mu(f^2) \le\frac{2}{\kappa}\int_\Omega|\nabla f|^2\,d\mu, \qquad f\in C^1(\overline\Omega).
\]

Consider \(r_t:=d\nu_t/d\mu=Z_\mu\rho_d(\cdot,t)e^\psi\). By \eqref{eq:closed_loop_wp_data_boundedness},
\[
    \frac{\operatorname*{ess\,sup}_{x\in\Omega}r_t(x)}
         {\operatorname*{ess\,inf}_{x\in\Omega}r_t(x)}
    \le\frac{L_d}{\ell_d}
       e^{\operatorname{osc}_{\overline\Omega}\psi},
\]
where \(\operatorname{osc}_{\overline\Omega}\psi :=\max_{\overline\Omega}\psi-\min_{\overline\Omega}\psi\). Set \(W_t:=-\ln r_t\). Then \(d\nu_t=e^{-W_t}\,d\mu\), and \(W_t\) is essentially bounded, with
\[
    \operatorname{ess\,osc}_{\Omega}W_t
    =\ln\frac{\operatorname*{ess\,sup}_{\Omega}r_t}
                  {\operatorname*{ess\,inf}_{\Omega}r_t}
    \le\ln\frac{L_d}{\ell_d}
       +\operatorname{osc}_{\overline\Omega}\psi.
\]
Since \(\mu\) satisfies the log-Sobolev inequality with coefficient \(2/\kappa\), the Holley--Stroock bounded perturbation principle~\cite[Theorem~A2]{Schlichting2019Mixtures} gives
\[
\begin{aligned}
    \operatorname{Ent}_{\nu_t}(f^2)
    &\le\frac{2}{\kappa}
       e^{\operatorname{ess\,osc}_{\Omega}W_t}
       \int_\Omega|\nabla f|^2\,d\nu_t\\
    &\le\frac{2L_d}{\kappa\ell_d}
       e^{\operatorname{osc}_{\overline\Omega}\psi}
       \int_\Omega|\nabla f|^2\,d\nu_t.
\end{aligned}
\]

To extend this inequality to \(f\in H^1(\Omega)\), choose \(f_j\in C^1(\overline\Omega)\) converging to \(f\) in \(H^1(\Omega)\). After passing to a subsequence, we also have \(f_j(x)\to f(x)\) for a.e. \(x\in\Omega\). Since \(\rho_d(\cdot,t)\le L_d\), the convergence \(f_j\to f\) in \(H^1(\Omega)\) implies
\[
\begin{aligned}
\|f_j-f\|_{L^2(\nu_t)}^2
&=\int_\Omega |f_j-f|^2\rho_d\,\mathrm{d}x\\
&\le L_d\|f_j-f\|_{L^2(\Omega)}^2\longrightarrow0,\\
\|\nabla f_j-\nabla f\|_{L^2(\nu_t)}^2
&=\int_\Omega |\nabla f_j-\nabla f|^2\rho_d\,\mathrm{d}x\\
&\le L_d\|\nabla f_j-\nabla f\|_{L^2(\Omega)}^2
\longrightarrow0.
\end{aligned}
\]
Hence \(\int_\Omega f_j^2\,d\nu_t\to\int_\Omega f^2\,d\nu_t\) and \(\int_\Omega|\nabla f_j|^2\,d\nu_t \to\int_\Omega|\nabla f|^2\,d\nu_t\). 

The convergence \(f_j\to f\) a.e. implies
\[
f_j^2\ln f_j^2+\frac1e
\longrightarrow f^2\ln f^2+\frac1e
\quad\text{a.e. in }\Omega.
\]
These functions are nonnegative because \(s\ln s\ge-1/e\). Fatou's lemma therefore gives
\[
\int_\Omega\left(f^2\ln f^2+\frac1e\right)d\nu_t
\le
\liminf_{j\to\infty}
\int_\Omega\left(f_j^2\ln f_j^2+\frac1e\right)d\nu_t.
\]
Since \(\nu_t(\Omega)=1\), subtracting \(1/e\) yields
\[
\int_\Omega f^2\ln f^2\,d\nu_t
\le\liminf_{j\to\infty}
\int_\Omega f_j^2\ln f_j^2\,d\nu_t.
\]
Moreover, writing \(m_j:=\int_\Omega f_j^2\,d\nu_t\) and \(m:=\int_\Omega f^2\,d\nu_t\), the convergence \(m_j\to m\) implies \(m_j\ln m_j\to m\ln m\). Hence,
\begin{equation}\label{eq:derive_process_of_log_Sobolev}
\begin{aligned}
\operatorname{Ent}_{\nu_t}(f^2)
&=\int_\Omega f^2\ln f^2\,d\nu_t-m\ln m\\
&\le\liminf_{j\to\infty}
\left(\int_\Omega f_j^2\ln f_j^2\,d\nu_t-m_j\ln m_j\right)\\
&=\liminf_{j\to\infty}\operatorname{Ent}_{\nu_t}(f_j^2)\\
& \le \liminf_{j\to \infty} \left(\frac{2L_d}{\kappa\ell_d} e^{\operatorname{osc}_{\overline\Omega}\psi}
       \int_\Omega|\nabla f_j|^2\,d\nu_t\right)\\
       & = \frac{2L_d}{\kappa\ell_d}e^{\operatorname{osc}_{\overline\Omega}\psi}\int_\Omega|\nabla f|^2\,d\nu_t
\end{aligned}
\end{equation}
for any \(f\in H^1(\Omega)\). Finally, taking \(f=\sqrt{\rho/\rho_d}\in H^1(\Omega)\) in \eqref{eq:derive_process_of_log_Sobolev}, then \(\int_\Omega f^2\,d\nu_t=1\), \(\operatorname{Ent}_{\nu_t}(f^2) =\operatorname{KL}\bigl(\rho\|\rho_d(\cdot,t)\bigr)\), and \(\int_\Omega|\nabla f|^2\,d\nu_t =\frac14\mathcal I\bigl(\rho\|\rho_d(\cdot,t)\bigr)\). Thus \eqref{eq:uniform_lsi_time_varying} holds, where \(C_{\mathrm{KL}}^{\mathrm{unif}} :=\frac{L_d}{2\kappa\ell_d}   e^{\operatorname{osc}_{\overline\Omega}\psi}\) is independent of \(t\) and \(\rho\).
\end{proof}

\begin{remark}
We can obtain an explicit upper bound for \(C_{\mathrm{KL}}^{\mathrm{unif}}\) as follows. Set \(\Delta_{\Omega} = \max_{\overline{\Omega}}\|x - x_0\|^2 - \min_{\overline{\Omega}}\|x - x_0\|^2 = \max_{\overline{\Omega}}\|x - x_0\|^2 \le d_{\Omega}^2\), then \(C_{\mathrm{KL}}^{\mathrm{unif}}(\kappa) =\frac{L_d}{2\kappa\ell_d}e^{\kappa\Delta_\Omega/2}\). Minimizing over \(\kappa\) gives \(\kappa=2/\Delta_\Omega\) and \(C_{\mathrm{KL}}^{\mathrm{unif}} =\frac{eL_d\Delta_\Omega}{4\ell_d}\). Moreover, note that
\[
    \frac{L_d}{\ell_d} = \frac{\max_m e^{a_m}}{\min_m e^{-a_m}} = e^{2\max_m a_m} = \left(\max_m e^{d_\Omega^2/(2\sigma_m^2)}\right)^2.
\]
Setting $\sigma_{\min} := \min_m \sigma_m$, then \({L_d}/{\ell_d} = \exp \left(d_{\Omega}^2 / \sigma_{\min}^2\right)\), and
\[
    C_{\mathrm{KL}}^{\mathrm{unif}} = \frac{\Delta_{\Omega}}{4}\exp \left(1 + \frac{d_{\Omega}^2}{\sigma_{\min}^2}\right) \le \frac{d_{\Omega}^2}{4}\exp \left(1 + \frac{d_{\Omega}^2}{\sigma_{\min}^2}\right).
\]
\end{remark}

\begin{theorem}[Exponential convergence]
\label{thm:time_varying_kl_convergence}
Under the assumptions of this section, let \(\rho_0\in L^2(\Omega)\) be a nonnegative probability density, and let \(\rho\) be the unique global weak solution given by Proposition~\ref{prop:closed_loop_wellposedness}. Then, for every \(t\ge0\),
\begin{equation}
    \begin{aligned}
    &\operatorname{KL}\bigl(\rho(\cdot,t)\|\rho_d(\cdot,t)\bigr)\\
    &\qquad\le e^{-Kt/C_{\mathrm{KL}}^{\mathrm{unif}}}
       \operatorname{KL}\bigl(\rho_0\|\rho_d(\cdot,0)\bigr),
    \end{aligned}
    \label{eq:time_varying_kl_exponential_decay}
\end{equation}
where \(C_{\mathrm{KL}}^{\mathrm{unif}}\) is the constant in Lemma~\ref{lem:uniform_lsi_time_varying}. In particular, the KL divergence converges to zero exponentially as \(t\to\infty\).
\end{theorem}

\begin{proof}
Fix \(T>0\), and set \(F(t):=\operatorname{KL}(\rho(t)\|\rho_d(t))\) and \(\mathcal I(t):=\mathcal I(\rho(t)\|\rho_d(t))\). Proposition~\ref{prop:KL_divergence_derivative} gives absolute continuity of \(F\) on \([0,T]\) and
\[
    F'(t)=-K\mathcal I(t),
    \quad \sqrt{\rho(t)/\rho_d(t)}\in H^1(\Omega)
    \quad\text{for a.e. }t\in(0,T).
\]
By Lemma~\ref{lem:uniform_lsi_time_varying},
\[
    F(t)\le C_{\mathrm{KL}}^{\mathrm{unif}}\mathcal I(t),
    \quad
    F'(t)\le-\frac{K}{C_{\mathrm{KL}}^{\mathrm{unif}}}F(t)
    \quad\text{for a.e. }t\in(0,T).
\]
The initial entropy is finite because
\(\rho_0\in L^2(\Omega)\) and \(0<\ell_d\le \rho_d(\cdot,0)\le L_d\). Gronwall's inequality therefore yields \eqref{eq:time_varying_kl_exponential_decay} on \([0,T]\). Since \(T>0\) is arbitrary and \(C_{\mathrm{KL}}^{\mathrm{unif}}\) is independent of \(T\), the estimate holds for all \(t\ge0\).
\end{proof}

\begin{remark}
The following results extend to general prescribed probability densities \(\rho_d\) as follows, with the domain and initial-data assumptions unchanged.
\begin{enumerate}
    \item In Proposition~\ref{prop:closed_loop_wellposedness}, the first two conclusions hold under \eqref{eq:closed_loop_wp_data_L_infty}--\eqref{eq:closed_loop_vd_data_L_infty}. The ratio regularity conclusion additionally requires \eqref{eq:reference_time_regularity}.

    \item Strict positivity in Proposition~\ref{prop:closed_loop_wellposedness} and KL dissipation in Proposition~\ref{prop:KL_divergence_derivative} hold under \eqref{eq:closed_loop_wp_data} and \eqref{eq:reference_transport_weak}. Theorem~\ref{thm:time_varying_kl_convergence} also applies if these conditions hold on every finite interval.

    \item Lemma~\ref{lem:uniform_lsi_time_varying} holds under
    \eqref{eq:closed_loop_wp_data_boundedness}.
\end{enumerate}
\end{remark}

\section[Estimated Feedback Law]{Estimated Feedback Law: Implementation and Conditional Robustness}
\label{section4}

The ideal feedback law \eqref{eq:diffusive_time_varying_control} uses the logarithmic gradient \(\nabla\ln\rho\) of the swarm density and a prescribed transport field \(v_d\). For implementation, we approximate these quantities using Gaussian kernel density estimation and filtered target velocities. Treating the estimator data as prescribed, we prove existence and uniqueness of a weak solution to the resulting Fokker--Planck system, and derive a conditional KL robustness estimate in terms of the feedback mismatch.

\subsection{Estimated Feedback Law Design}

\paragraph{Estimation of the Logarithmic Density Gradient}

For fixed \(N\ge1\) and bandwidth \(h>0\), the full-space Gaussian
kernel density estimate (KDE) ~\cite{Rosenblatt1956} is
\begin{equation}
    \widetilde\rho_h^N(x,t)=\frac1N\sum_{i=1}^N K_h(x-x_i(t)),
    \quad
    K_h(z)=\frac{1}{2\pi h^2}e^{-|z|^2/(2h^2)}.\label{eq:kde_estimate}
\end{equation}
Its normalization on the task domain is
\begin{equation}
    \hat\rho_h^N(x,t):=\frac{\widetilde\rho_h^N(x,t)}{C_h^N(t)},\quad C_h^N(t):=\int_\Omega\widetilde\rho_h^N(z,t)\,dz. \label{eq:normalized_kde_estimate}
\end{equation}
The Gaussian kernel is strictly positive, so both \(C_h^N(t)\) and
\(\hat\rho_h^N\) are positive. Since the normalizing factor is
independent of \(x\),
\begin{equation}
    \nabla\ln\hat\rho_h^N=\nabla\ln\widetilde\rho_h^N.
    \label{eq:kde_normalization_score_invariance}
\end{equation}
Define \(a_i(x,t):=\frac{K_h(x-x_i(t))}{\sum_jK_h(x-x_j(t))}\), then \(a_i\ge0\), \(\sum_i a_i=1\), and differentiation gives \(\nabla\ln\hat\rho_h^N(x,t) =-\frac1{h^2}\left(x-\sum_i a_i(x,t)x_i(t)\right)\). For \(x,x_i(t)\in\overline\Omega\), convexity implies that the \(\sum_i a_i x_i\) also lies in \(\overline\Omega\). Since \(|x-x_i(t)|\le d_\Omega\), the triangle inequality gives
\begin{equation}
    \|\nabla\ln\hat\rho_h^N(\cdot,t)\|_{L^\infty(\Omega)}
    \le\frac{d_\Omega}{h^2}.
    \label{eq:kde_score_boundedness}
\end{equation}

\paragraph{Estimation of the Reference Transport Velocity}

The weighted combination of component velocities has been used in previous studies \cite{lipman2023flow,rapakoulias2025steering}. We adopt this established construction to obtain an estimation of the reference transport velocity \(v_d\). Consider first the full space Gaussian mixture \(\rho_d^{\infty}(x,t) = \sum_m w_m \phi_m^{\infty}(x,t)\), where \(\phi_m^\infty(x,t):=(2\pi\sigma_m^2)^{-1} \exp\left(-|x-y_m(t)|^2/(2\sigma_m^2)\right)\), and \(w_m\) are given in \eqref{eq:task_demand_density}. The chain rule gives \(\partial_t\phi_m^\infty =-\nabla\cdot(\phi_m^\infty \ \dot{y}_m)\). Multiplying the component equations by \(w_m\) and summing gives
\[
    \partial_t \rho_d^{\infty}
    =-\nabla\cdot\left(\sum_m w_m \ \phi_m^{\infty} \ \dot{y}_m\right).
\]
Since \(\rho_d^{\infty}(x,t)>0\), define \(c(x,t):=\sum_m\frac{w_m \phi_m^{\infty}(x,t)}{\rho_d^{\infty}(x,t)}\dot{y}_m(t)\), then this field satisfies \(\partial_t \rho_d^{\infty}+\nabla\cdot(\rho_d^{\infty}c)=0\) and therefore transports the mixture density \(\rho_d^{\infty}\). 

For the truncated components on the task domain, we use the same density-weighted construction to define a kinematic approximation of \(v_d\)
\begin{equation}
    v_d^{\mathrm{kin}}(x,t):=\sum_m\gamma_m(x,t)\dot y_m(t).
    \label{eq:vd_estimate}
\end{equation}
where \(\gamma_m(x,t):=\frac{w_m\phi_m(x,t)}{\rho_d(x,t)}\) satisfy \(\gamma_m\ge0\) and \(\sum_m\gamma_m=1\). 

For implementation, we estimate \(\dot y_m(t)\) in \eqref{eq:vd_estimate} from measured target positions \(y_m\) using a first-order filtered differentiator. For each target \(m\), introduce an auxiliary filtered position \(z_m\in\mathbb R^2\) satisfying
\begin{equation}
    \dot z_m=\lambda_m(y_m-z_m), \qquad z_m(0)=y_m(0), \label{eq:target_position_filter}
\end{equation}
where \(\lambda_m>0\) is the filter bandwidth. The target-velocity estimate is defined by
\begin{equation}
    \widehat{\dot y}_m:=\lambda_m(y_m-z_m)=\dot z_m.
    \label{eq:filtered_target_velocity}
\end{equation}
Then the implemented transport field is
\begin{equation}
    \hat v_d(x,t) =\sum_m\gamma_m(x,t)\widehat{\dot y}_m(t).
    \label{eq:filtered_vd_estimate}
\end{equation}
And the estimated feedback is
\begin{equation}
    \boxed{\begin{aligned} \hat u={}&\hat v_d+K\nabla\ln\rho_d - (K-D)\nabla\ln\hat\rho_h^N.
    \end{aligned}}
    \label{eq:estimated_feedback_control}
\end{equation}

\subsection[Weak Well-Posedness]{Weak Well-Posedness for a Prescribed Estimated Drift}
\label{subsec:estimated_wellposedness}

\begin{proposition}
\label{prop:prescribed_drift_wellposedness}
Let \(\Omega\) be the task domain of Section~\ref{section2},
\(T>0\), and \(D>0\). If \(u\in L^\infty(\Omega\times(0,T);\mathbb R^2)\), \(\rho_0\in L^2(\Omega)\), \(\rho_0\ge 0\) and \(\int_\Omega\rho_0\,dx=1\), then \eqref{eq:rho_u_FP} has a unique weak solution
\begin{equation}
    \begin{gathered}
    \rho\in C([0,T];L^2(\Omega))\cap L^2(0,T;H^1(\Omega)),\\
    \partial_t\rho\in L^2(0,T;(H^1(\Omega))^*),
    \qquad \rho(0)=\rho_0,
    \end{gathered}
    \label{eq:prescribed_drift_solution_class}
\end{equation}
satisfying, for every \(\psi\in H^1(\Omega)\) and a.e. \(t\),
\begin{equation}
    \langle\partial_t\rho,\psi\rangle
    =-D\int_\Omega\nabla\rho\cdot\nabla\psi\,dx
      +\int_\Omega\rho u\cdot\nabla\psi\,dx.
    \label{eq:prescribed_drift_weak_form}
\end{equation}
For every \(t\in[0,T]\), \(\rho(t)\ge0\) a.e. and
\(\int_\Omega\rho(t)\,dx=1\). For two initial data with the same drift \(u\), their solutions satisfy
\begin{equation}
    \|\rho_1(t)-\rho_2(t)\|_2^2 \le e^{\frac{U_T^2t}{D} }\|\rho_{1,0}-\rho_{2,0}\|_2^2, \quad U_T:=\|u\|_{L^\infty(\Omega\times(0,T))}.
    \label{eq:prescribed_drift_continuous_dependence}
\end{equation}
If the drift is bounded on every finite time interval, these
solutions define a unique global weak solution.
\end{proposition}

\begin{proof}
Set \(V:=H^1(\Omega)\), \(H:=L^2(\Omega)\), and
\[
    a_t(f,g):=D\int_\Omega\nabla f\cdot\nabla g\,dx
       -\int_\Omega fu\cdot\nabla g\,dx.
\]
The form is measurable in \(t\) and satisfies
\[
    \begin{aligned}
    |a_t(f,g)|&\le(D+U_T)\|f\|_V\|g\|_V,\\
    a_t(f,f)&\ge\frac D2\|\nabla f\|_2^2
                    -\frac{U_T^2}{2D}\|f\|_2^2.
    \end{aligned}
\]
Adding \((U_T^2/(2D)+D/2)\|f\|_2^2\) to the second inequality gives 
\[
    a_t(f,f) + \left( \frac{U_T^2}{2D} + \frac{D}{2}\right)\|f\|_2^2 \ge \frac{D}{2}\|f\|_V^2
\]
Hence $a_t$ is a non-autonomous closed form on $H$. Lions' variational theorem, as used in
Proposition~\ref{prop:closed_loop_wellposedness}, therefore yields a
unique solution in \(L^2(0,T;V)\cap H^1(0,T;V^*)\), with initial
value \(\rho_0\) in \(H\). The embedding into \(C([0,T];H)\)
gives \eqref{eq:prescribed_drift_solution_class}, and the definition
of the form gives \eqref{eq:prescribed_drift_weak_form}.

For \(\rho^-:=\max\{-\rho,0\}\), the truncation chain rule and
the test function \(-\rho^-\) give
\[
    \frac12\frac d{dt}\|\rho^-\|_2^2
       +D\|\nabla\rho^-\|_2^2
    \le\frac D2\|\nabla\rho^-\|_2^2
       +\frac{U_T^2}{2D}\|\rho^-\|_2^2.
\]
Gronwall's inequality and \(\rho_0^-=0\) imply nonnegativity.
Testing with \(1\) gives conservation of mass. These conclusions
hold at every time for the continuous \(H\)-representative.
The same energy calculation for \(\rho_1-\rho_2\) proves
\eqref{eq:prescribed_drift_continuous_dependence}. Finally,
uniqueness makes solutions on overlapping finite intervals agree.
\end{proof}

For \eqref{eq:estimated_feedback_control}, fix \(N\ge1\) and \(h>0\), prescribe continuous paths \(x_i:[0,\infty)\to\overline\Omega\), and retain the target-trajectory assumptions of Section~\ref{section3}. Then the estimates \(\widehat{\dot y}_m(t)\) are continuous and hence bounded on \([0,T]\). The reference assumptions give \(\|\nabla\ln\rho_d\|_{L^\infty(\Omega\times(0,T))} \le\frac{1}{\ell_d} \|\nabla\rho_d\|_{L^\infty(\Omega\times(0,T))} < \infty\), while \eqref{eq:kde_score_boundedness} bounds \(\nabla\ln\hat\rho_h^N\). Using \(\gamma_m\ge0\) and \(\sum_m\gamma_m=1\), we obtain
\[
    \begin{aligned}
    \|\hat u\|_{L^\infty(\Omega\times(0,T))}
    &\le\max_m\sup_{0\le t\le T}|\widehat{\dot y}_m(t)| + K\|\nabla\ln\rho_d\|_{L^\infty(\Omega\times(0,T))}\\
    &\quad+|K-D|\frac{d_\Omega}{h^2}<\infty.
    \end{aligned}
\]
Thus Proposition~\ref{prop:prescribed_drift_wellposedness} applies with \(u=\hat u\). Since \(T>0\) is arbitrary, uniqueness yields a unique global weak solution.

\subsection[Conditional KL Robustness]{Conditional Robustness of the KL Convergence}
\label{subsec:conditional_robustness}

We now study the robustness of KL convergence with respect to estimation errors. Let \(\rho\) be the global weak solution corresponding to \(u = \hat u\) in \eqref{eq:estimated_feedback_control}, with initial density \(\rho_0\) satisfying Proposition~\ref{prop:prescribed_drift_wellposedness}.  Since \(\rho\) may vanish, for a.e. \(t>0\) we define
\[
    s_\rho:=
    \begin{cases}
       \nabla\rho/\rho,&\rho>0,\\
       0,&\rho=0.
    \end{cases}
\]
The identity \(\nabla\rho=0\) a.e. on \(\{\rho=0\}\) then implies \(\rho s_\rho=\nabla\rho\) a.e. Define the logarithmic gradient estimation error \(e_\rho\), the transport-velocity error \(e_v\), and the combined disturbance \(d\) by
\begin{equation}
    \begin{aligned}
    e_\rho&:=\nabla\ln\hat\rho_h^N-s_\rho,\qquad e_v:=\hat v_d-v_d,\\
    d&:=e_v-(K-D)e_\rho\\
     &=\hat u- \left(v_d+K\nabla\ln\rho_d - (K-D)s_\rho \right).
    \end{aligned}\label{eq:combined_estimation_disturbance}
\end{equation}

\begin{theorem}[Conditional robust KL bound]
\label{thm:robust_kl_bound}
Let \(D>0\), \(K>0\), and let \(\rho\) be the global weak solution \(\rho\) introduced above. Assume that the prescribed \((\rho_d,v_d)\) satisfies \eqref{eq:closed_loop_wp_data} and \eqref{eq:reference_transport_weak} on every finite time interval, with the same bounds \(\ell_d,L_d\) for all time. Let \(C_{\mathrm{KL}}^{\mathrm{unif}}\) be the constant provided by Lemma~\ref{lem:uniform_lsi_time_varying}. If the disturbance in \eqref{eq:combined_estimation_disturbance} satisfies
\begin{equation}
    M(t):=\|d(\cdot,t)\|_{\rho(t)}^2,
    \qquad \int_0^T M(t)\,dt<\infty
    \quad\text{for every }T>0.
    \label{eq:conditional_disturbance_integrability}
\end{equation}
Then the following properties hold:
\begin{enumerate}
    \item For every \(T>0\), \(\sqrt{\rho/\rho_d}\in L^2(0,T;H^1(\Omega))\), and \(F(t):=\operatorname{KL}(\rho(t)\|\rho_d(t))\) is absolutely continuous on \([0,T]\). With \(q:=\rho/\rho_d\), it obeys
\begin{equation}
    F'(t)=-K\mathcal I(t)    +2\int_\Omega\rho_d\sqrt q\,d\cdot\nabla\sqrt q\,dx \quad\text{for a.e. }t, \label{eq:perturbed_kl_identity}
\end{equation}
where \(\mathcal I(t)\) is defined by
\eqref{eq:relative_fisher_information}. For positive regular
densities, the second term equals
\(\int_\Omega\rho d\cdot\nabla\ln(\rho/\rho_d)\,dx\).

    \item With \(\lambda:=K/(2C_{\mathrm{KL}}^{\mathrm{unif}})\), for every \(t \ge 0\),
\begin{equation}
    F(t)\le e^{-\lambda t}F(0)    +\frac1{2K}\int_0^t e^{-\lambda(t-s)}M(s)\,ds. \label{eq:robust_kl_convolution_bound}
\end{equation}
    \item If, additionally, there exist constants \(\varepsilon_v,\varepsilon_\rho\ge0\) such that
\begin{equation}
    \|e_v(\cdot,t)\|_{\rho(t)}\le\varepsilon_v, \quad \|e_\rho(\cdot,t)\|_{\rho(t)}\le\varepsilon_\rho \quad\text{a.e. }t\ge0, \label{eq:uniform_estimation_error_bounds}
    \end{equation}
    then, with \(\varepsilon:=\varepsilon_v+|K-D|\varepsilon_\rho\),
\begin{equation}
    F(t)\le e^{-\lambda t}F(0)    +\frac{C_{\mathrm{KL}}^{\mathrm{unif}}\varepsilon^2}{K^2}       (1-e^{-\lambda t}), \label{eq:robust_kl_ultimate_bound}
\end{equation}
and \(\limsup_{t\to\infty}F(t)
\le C_{\mathrm{KL}}^{\mathrm{unif}}\varepsilon^2/K^2\).
\end{enumerate} 
\end{theorem}

\begin{proof}
    1. Fix \(T>0\). Set \(q:=\rho/\rho_d\). The same argument as in the proof of Proposition~\ref{prop:closed_loop_wellposedness} gives \(q\in C([0,T];L^2(\Omega))\cap L^2(0,T;H^1(\Omega))\) and \(\partial_tq\in L^2(0,T;(H^1(\Omega))^*)\). Since \(\rho s_\rho=\nabla\rho\) a.e. and \(\rho_d\nabla q =\nabla\rho-\rho\nabla\ln\rho_d\), we obtain \(\hat J:=\rho\hat u-D\nabla\rho =\rho v_d-K\rho_d\nabla q+\rho d \) a.e. in  \(\Omega\times(0,T)\).

Fix \(\varepsilon \in (0,1]\), set \(f_\varepsilon:=\rho+\varepsilon\rho_d\), \(q_\varepsilon:=q+\varepsilon\), \(\eta(r):=r\ln r-r+1\) and \(E_\varepsilon:=\int_\Omega\rho_d\eta(q_\varepsilon)\,dx\). As in the proof of Proposition~\ref{prop:KL_divergence_derivative}, we write
\[
    E_\varepsilon =\int_\Omega    \bigl[\eta(f_\varepsilon)  -f_\varepsilon\ln\rho_d+\rho_d-1\bigr]\,dx.
\]
Since \(f_\varepsilon\ge\varepsilon\ell_d>0\), the same variational chain-rule and time product-rule argument as in the proof of Proposition~\ref{prop:KL_divergence_derivative} shows that \(E_\varepsilon\) is absolutely continuous on \([0,T]\), with
\[
    E_\varepsilon' =\langle\partial_t f_\varepsilon,\ln(q+\varepsilon)\rangle    +\int_\Omega(1-q-\varepsilon)\partial_t\rho_d\,dx.
\]
Here, \(\langle\partial_t f_\varepsilon,\psi\rangle = \int_\Omega J_\varepsilon\cdot\nabla\psi\,dx\ (\psi\in H^1(\Omega))\), so \(\langle\partial_t f_\varepsilon,\ln(q+\varepsilon)\rangle =\int_\Omega J_\varepsilon\cdot    \frac{\nabla q}{q+\varepsilon}\,dx\), where \(J_\varepsilon :=\hat J+\varepsilon\rho_dv_d = \rho_d(q+\varepsilon)v_d -K\rho_d\nabla q+\rho\,d\). Testing \eqref{eq:reference_transport_weak} with \(1-q-\varepsilon\) gives \(\int_\Omega(1-q-\varepsilon)\partial_t\rho_d\,dx = -\int_\Omega\rho_dv_d\cdot\nabla q\,dx\). Thus
\begin{equation}
    E_\varepsilon' =-K\int_\Omega\rho_d \frac{|\nabla q|^2}{q+\varepsilon}\,dx +\int_\Omega\rho\,d\cdot \frac{\nabla q}{q+\varepsilon}\,dx.
    \label{eq:regularized_perturbed_entropy}
\end{equation}
Note that \(\int_0^T\!\int_\Omega \rho_d\frac{|\nabla q|^2}{q+\varepsilon}\,dx\,dt \le\frac{L_d}{\varepsilon} \int_0^T\|\nabla q(\cdot,t)\|_2^2\,dt<\infty\), \(\left|\int_\Omega\rho\,d\cdot \frac{\nabla q}{q+\varepsilon}\,dx\right| \le\sqrt{\frac{L_d}{\varepsilon}}\, \|d(\cdot,t)\|_{\rho(t)}\|\nabla q(\cdot,t)\|_2\). Thus each integral in \eqref{eq:regularized_perturbed_entropy} is integrable in time, and by the absolute continuity of \(E_{\varepsilon}\), integrating over \(0\le s\le t \le T\) yields
\begin{equation}
    E_{\varepsilon}(t) - E_{\varepsilon}(s) = -K \int_s^t A_{\varepsilon}(\tau)\, d \tau + \int_s^t \int_{\Omega}\rho d \cdot \frac{\nabla q}{q + \varepsilon}\, dx \, d\tau\label{eq:regularized_perturbed_entropy_integral_form}
\end{equation}
for every \(0\le s\le t \le T\), where \(A_\varepsilon(t):=\int_\Omega \frac{\rho_d|\nabla q|^2}{q+\varepsilon}\,dx\).
 
2. By Young's inequality,
\[
    \begin{aligned}
    \left|\int_\Omega\rho d\cdot
                    \frac{\nabla q}{q+\varepsilon}\,dx\right|
    &\le\frac K2 A_\varepsilon+\frac1{2K}\int_\Omega
          \rho_d\frac{q^2}{q+\varepsilon}|d|^2\,dx\\
    &\le\frac K2 A_\varepsilon+\frac1{2K}M.
    \end{aligned}
\]
Combining this estimate with \eqref{eq:regularized_perturbed_entropy} gives \(E_\varepsilon'+(K/2)A_\varepsilon\le M/(2K)\). Integrating over \([s,t]\) yields
\begin{equation}
    E_\varepsilon(t)+\frac K2\int_s^t A_\varepsilon(\tau)\,d\tau
    \le E_\varepsilon(s)+\frac1{2K}\int_s^t M(\tau)\,d\tau.
    \label{eq:regularized_entropy_bound}
\end{equation}
As in the proof of Proposition~\ref{prop:KL_divergence_derivative}, since \(0\le \rho_d \eta(q+\varepsilon) \le 2L_d(q^2 + 1)\) for any \(\varepsilon \in (0,1]\), dominated convergence gives, for every \(t \in [0,T]\),
\[
    E_\varepsilon(t)\to \int_\Omega\rho_d\eta(q)\,dx=F(t).
\]
The same bound at \(t = 0\) gives \(E_\varepsilon(0) \le 2L_d\left(\|q(0)\|_{L^2(\Omega)}^2+|\Omega|\right)<\infty\) for any \(\varepsilon\in (0,1]\). Taking \(s = 0,t = T\) in \eqref{eq:regularized_entropy_bound}, and using
\(E_\varepsilon(T)\ge0\), yields
\[
\begin{gathered}
    \int_0^T A_\varepsilon(\tau)\,d\tau \le \frac2K E_\varepsilon(0) +\frac1{K^2}\int_0^T M(\tau)\,d\tau\\
 \le \frac{4L_d}{K} \left(\|q(0)\|_{L^2(\Omega)}^2+|\Omega|\right) + \frac1{K^2}\int_0^T M(\tau)\,d\tau  < \infty.
\end{gathered}
\]
Since \(\rho_d\frac{|\nabla q|^2}{q+\varepsilon}\ge 0\) in \(A_{\varepsilon}\) increases as \(\varepsilon\downarrow 0\), the monotone convergence then yields
\[
    \int_0^T A_{\varepsilon}(t) \, dt \to \int_0^T\!\int_\Omega\rho_d\frac{|\nabla q|^2}{q}\,dx\,dt
    <\infty,\quad \varepsilon\downarrow 0,
\]
where the quotient is zero on \(\{q=0\}\). Moreover, the same weak compactness argument as in the proof of Proposition~\ref{prop:KL_divergence_derivative} gives \(\sqrt q\in L^2(0,T;H^1(\Omega))\). The Sobolev chain rule then gives \(\rho_d\frac{|\nabla q|^2}{q} =4\rho_d|\nabla\sqrt q|^2\) a.e., and
\[
    \mathcal I(t) =4\int_\Omega\rho_d|\nabla\sqrt q|^2\,dx =\int_\Omega\rho_d\frac{|\nabla q|^2}{q}\,dx \quad\text{for a.e. }t.
\]
Hence \(\mathcal{I} \in L^1(0,T)\). Applying monotone convergence on each subinterval gives, as \(\varepsilon \downarrow 0\),
\[
    \int_s^t A_{\varepsilon}(\tau) \, d\tau \to \int_s^t\mathcal{I}(\tau)\,d\tau,\quad 0\le s\le t \le T.
\]

For the disturbance term, note that 
\[
    \rho d\cdot\frac{\nabla q}{q+\varepsilon} =2\rho_d\sqrt q\,d\cdot\nabla\sqrt q\, \frac{q}{q+\varepsilon}.
\]
This converges a.e. to \(2\rho_d\sqrt q\,d\cdot\nabla\sqrt q\) as \(\varepsilon \downarrow 0\), and its absolute value is bounded by \(G :=2\rho_d\sqrt q\,|d|\,|\nabla\sqrt q|\) such that \(\int_0^T\!\int_\Omega G\,dx\,dt \le \left(\int_0^T M(t)\,dt\right)^{1/2} \left(\int_0^T\mathcal I(t)\,dt\right)^{1/2} <\infty\). Thus dominated convergence yields, for \(0\le s\le t\le T\), as \(\varepsilon\downarrow 0\),
\[
    \int_s^t\!\int_\Omega \rho d\cdot\frac{\nabla q}{q+\varepsilon}\,dx\,d\tau \to 2\int_s^t\!\int_\Omega \rho_d\sqrt q\,d\cdot\nabla\sqrt q\,dx\,d\tau.
\]
Passing to the limit in \eqref{eq:regularized_perturbed_entropy_integral_form} gives, for every \(0\le s\le t\le T\),
\[
    F(t)+K\int_s^t\mathcal I(\tau)\,d\tau =F(s)+2\int_s^t\!\int_\Omega \rho_d\sqrt q\,d\cdot\nabla\sqrt q\,dx\,d\tau.
\]
Since both time integrands belong to \(L^1(0,T)\), \(F\) is absolutely continuous on \([0,T]\), and \eqref{eq:perturbed_kl_identity} follows.

3. Finally, since \(\sqrt q\in L^2(0,T;H^1(\Omega))\), by Lemma~\ref{lem:uniform_lsi_time_varying}, \eqref{eq:uniform_lsi_time_varying} holds. Cauchy--Schwarz and Young's inequality give \(\left|2\int_\Omega\rho_d\sqrt q\,d\cdot\nabla\sqrt q\,dx\right| \le\sqrt{M(t)}\sqrt{\mathcal I(t)} \le\frac1{2K}M(t)+\frac K2\mathcal I(t)\). So 
\begin{align*}
F'(t) &=-K\mathcal I(t)
  +2\int_\Omega\rho_d\sqrt q\,d\cdot\nabla\sqrt q\,dx\\
&\le -\frac K2\mathcal I(t)+\frac1{2K}M(t).
\end{align*}
Combining this with \eqref{eq:uniform_lsi_time_varying} yields 
\[
    F'(t)+\lambda F(t)\le\frac1{2K}M(t)
    \quad\text{for a.e. }t\in(0,T).
\]
where \(\lambda = K / (2C_{\mathrm{KL}}^{\mathrm{unif}})\). The integrating-factor argument proves
\eqref{eq:robust_kl_convolution_bound}. Under \eqref{eq:uniform_estimation_error_bounds}, the triangle inequality gives \(M(t)\le(\varepsilon_v+|K-D|\varepsilon_\rho)^2\) a.e., which proves \eqref{eq:robust_kl_ultimate_bound}. Since \(T\) is arbitrary, the conclusions hold globally.
\end{proof}

When \(d=0\), \eqref{eq:perturbed_kl_identity} recovers \eqref{eq:time_varying_kl_dissipation}. In practice, the control magnitude is typically bounded. For a prescribed maximum admissible speed \(u_{\max}>0\), we enforce \(|u_{\mathrm a}|\le u_{\max}\) by defining 
\begin{equation}
    u_{\mathrm a}:= \begin{cases}
        \hat u,&|\hat u|\le u_{\max},\\
        u_{\max}\dfrac{\hat u}{|\hat u|},
        &|\hat u|>u_{\max}.
    \end{cases}. \label{eq:ua}
\end{equation}
Since \(u_{\mathrm a}\) is measurable and \(|u_{\mathrm a}|\le u_{\max}\), Proposition~\ref{prop:prescribed_drift_wellposedness} yields a unique global weak solution \(\rho\) with drift \(u_{\mathrm a}\) and the same initial density \(\rho_0\). Define \(e_\rho\), \(e_v\), and \(d\) as above using this solution. The total disturbance then includes the saturation error \(d_{\mathrm{tot}}:=d+(u_{\mathrm a}-\hat u)\). The entropy identity and convolution bound apply with \(d_{\mathrm{tot}}\) replacing \(d\), provided \(\int_0^T\|d_{\mathrm{tot}}\|_{\rho}^2\,dt<\infty\) for every \(T>0\). A uniform bound \(\|d_{\mathrm{tot}}\|_\rho\le\varepsilon_{\mathrm{tot}}\) yields the ultimate estimate with \(\varepsilon_{\mathrm{tot}}\) replacing \(\varepsilon\).

\section{Simulation}\label{sec:simulation}

In this section, we evaluate the proposed controller in a moving-target tracking experiment, compare its tracking performance with and without \(\hat v_d\) over different random seeds, and then compare it with the density-tracking feedback of Zheng~\cite{Zheng2021FieldEstimation} on the same task.

\subsection{Experimental setup}
\label{subsec:experimental_setup}

We simulate \(N=240\) agents in \(\overline\Omega=[0,20]\times[0,10]\) for \(T=24\) seconds with time step \(\Delta t=0.03\). Initial positions are sampled from a Gaussian distribution with mean \((0.7,5.0)\) and componentwise standard deviations \((0.25,0.65)\), then reflected into the task region. We use \(K=0.9\), \(D=9.8\times10^{-5}\), KDE bandwidth \(h=0.30\), and filter bandwidth \(\lambda_m=2\), with \(z_m(0)=y_m(0)\). The target components have weights \(w_m=1/3\) and spatial scales \((\sigma_1,\sigma_2,\sigma_3)=(0.70,0.75,0.70)\).

For \(s=t/T\), define \(q_1(t)=4.5+12[s+0.1\sin(\pi s)]\), \(q_2(t)=4.5+12[s-0.1\sin(\pi s)]\), and \(q_3(t)=4.5+12s\). The target trajectories are
\begin{equation}
    \begin{aligned}
    y_1(t)&=\left(q_1(t),
      5-3\cos\frac{\pi(q_1(t)-4.5)}{12}\right),\\
    y_2(t)&=\left(q_2(t),
      5+3\cos\frac{\pi(q_2(t)-4.5)}{12}\right),\\
    y_3(t)&=\left(q_3(t),5\right).
    \end{aligned}
    \label{eq:complex_target_paths}
\end{equation}
The desired density \(\rho_d\) is given by \eqref{eq:task_demand_density}, with components centered at the target positions given in \eqref{eq:complex_target_paths}.

To assess whether the static-density feedback can track moving targets, we compare \eqref{eq:estimated_feedback_control} with
\begin{equation}
    \hat u_0 =K\nabla\ln\rho_d-(K-D)\nabla\ln\hat\rho_h^N,\label{eq:hat_u0}
\end{equation}
which is the proposed controller without \(\hat{v}_d\), corresponding to its static-density form. For each controller, let \(u_a\) denote the applied drift using the saturation rule \eqref{eq:ua} with maximum admissible speed \(u_{\max}=3\). The agent positions are then updated using the reflected Euler-Maruyama scheme:
\begin{equation}
    x_i^{k+1}=\mathcal R_\Omega\!\left[  x_i^k+u_{\mathrm a}(x_i^k,t_k)\Delta t  +\sqrt{2D\Delta t}\,\omega_i^k\right],\label{eq:simulation_em}
\end{equation}
where \(\omega_i^k\sim\mathcal N(0,I_2)\), and \(x_i^k\) is the position of agent \(i\) at \(t_k = k\Delta t\). The coordinatewise reflection map \(\mathcal R_\Omega:\mathbb R^2\to\overline\Omega\) is defined by \(\mathcal R_\Omega(q) =\bigl(\mathcal R_{[0,20]}(q_1), \mathcal R_{[0,10]}(q_2)\bigr)^\top\), where \(\mathcal R_{[a,b]}(s)\) is obtained by initializing \(r=s\) and applying \(r\leftarrow \begin{cases}
2a-r,&r<a,\\
2b-r,&r>b,
\end{cases}\) until \(r\in[a,b]\). We repeat the experiment using random seeds 23--32. Within each run, all controllers use the same initial agent positions and Gaussian random vectors. Across runs, these random quantities are regenerated, while all other simulation settings remain fixed. 

Density errors are evaluated every \(10 \Delta t = 0.30\) seconds on a uniform \(88\times44\) grid. With grid spacings \(\Delta x, \Delta y\) and \(\Delta A = \Delta x \Delta y\), let \(\hat\rho_{pq}(t)\) and \(\rho_{d,pq}(t)\) denote the KDE and desired-density values at grid point \((p,q)\), respectively. We normalize both grid densities as \(\bar\rho_{pq}(t)= \frac{\hat\rho_{pq}(t)} {\sum_{r,s}\hat\rho_{rs}(t)\Delta A}\), \(\bar\rho_{d,pq}(t)=\frac{\rho_{d,pq}(t)} {\sum_{r,s}\rho_{d,rs}(t)\Delta A}\). The discrete KL error is then defined by
\begin{equation}
    F_{\mathrm{num}}(t)
    :=\sum_{p,q}\bar\rho_{pq}(t)
       \ln\frac{\bar\rho_{pq}(t)}
               {\bar\rho_{d,pq}(t)}\,\Delta A.
    \label{eq:discrete_kl_metric}
\end{equation}
We report the pointwise mean \(\overline F_{\mathrm{num}}\) across the ten runs. For the late-time error, we first average the samples in \(t\in[16.8,24]\), the final \(30\%\) of the simulation, within each run, and then compute the mean and sample standard deviation across runs.

\subsection{Tracking performance}
\label{subsec:simulation_complex_targets}

\begin{figure*}[!t]
    \centering
    \captionsetup[subfigure]{skip=3pt}
    \begin{subfigure}{\textwidth}
        \centering
        \includegraphics[width=\linewidth]{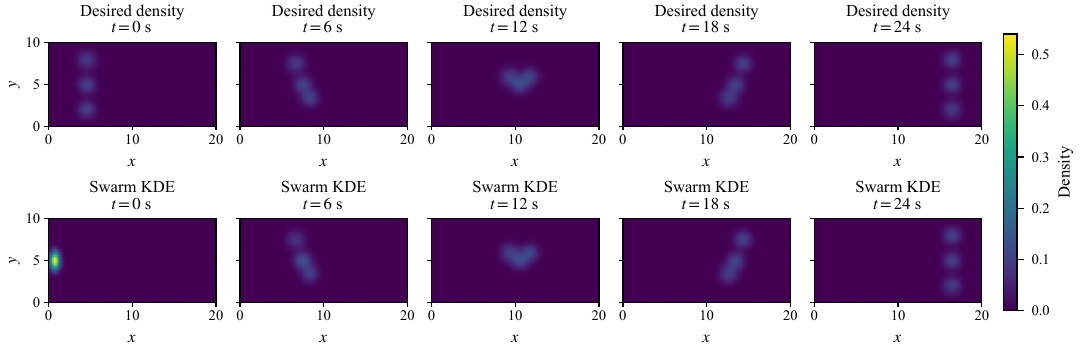}
        \caption{Desired-density and swarm-KDE evolution through the crossing interval.}
        \label{fig:complex_density_panel}
    \end{subfigure}\\[2mm]
    \begin{subfigure}{\textwidth}
        \centering
        \includegraphics[width=\linewidth]{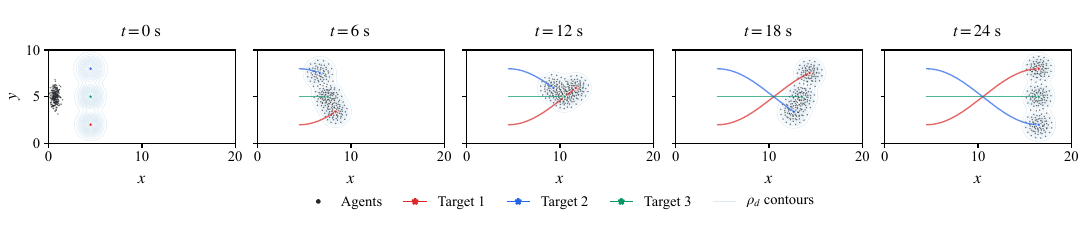}
        \caption{Corresponding agent positions and target paths.}
        \label{fig:complex_particles_panel}
    \end{subfigure}
    \caption{Representative experiment under the proposed controller, using seed 23, at \(t=0,6,12,18,24\) seconds. In (b), dots denote agents, colored stars denote targets, colored curves show their traversed paths, and blue contours indicate the desired task-demand density.}
    \label{fig:complex_swarm_snapshots}
\end{figure*}

Figure~\ref{fig:complex_swarm_snapshots}\subref{fig:complex_density_panel} shows the desired density and the swarm KDE at five representative times for seed 23. Figure~\ref{fig:complex_swarm_snapshots}\subref{fig:complex_particles_panel} shows the corresponding agent positions and target paths, together with contours of the desired density. The swarm redistributes as the targets approach, cross, and separate, suggesting that the swarm is able to track multiple moving targets without explicit agent--target assignments. Supplementary Video~1 shows the corresponding continuous evolution for \(h=0.30\).

\begin{figure*}[!t]
    \centering
    \captionsetup[subfigure]{skip=3pt}
    \begin{subfigure}{0.45\textwidth}
        \centering
        \includegraphics[width=\linewidth]{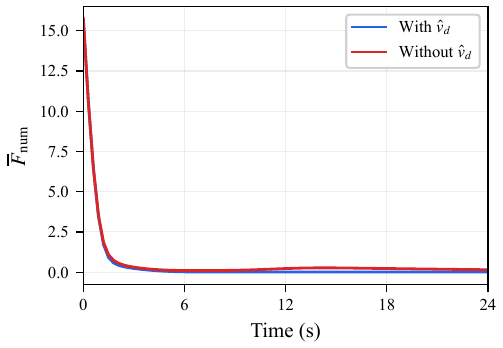}
        \caption{Mean KL error.}
        \label{fig:complex_kl}
    \end{subfigure}\hfill
    \begin{subfigure}{0.45\textwidth}
        \centering
        \includegraphics[width=\linewidth]{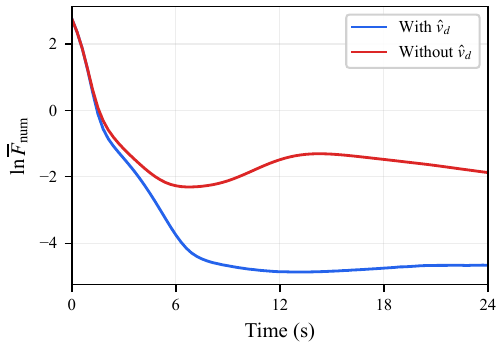}
        \caption{Logarithm of the mean KL error.}
        \label{fig:complex_log_kl}
    \end{subfigure}
    \caption{Tracking errors over random seeds 23--32. Panel (a) shows the pointwise mean KL error across the ten runs, and panel (b) shows its logarithm. Each controller pair uses the same initial particles, Brownian increments, and speed limit.}
    \label{fig:complex_kl_comparison}
\end{figure*}

Figure~\ref{fig:complex_kl_comparison} compares the mean tracking errors of controllers \eqref{eq:estimated_feedback_control} and \eqref{eq:hat_u0} over the ten random seeds. 
Figure~\ref{fig:complex_kl_comparison}\subref{fig:complex_kl} shows that, under the proposed controller, the mean KDE-based tracking error decreases rapidly during the initial transient and subsequently remains small. This behavior is qualitatively consistent with the convergence and conditional robustness analysis in Sections~\ref{section3} and~\ref{section4}. For a clear comparison, Figure~\ref{fig:complex_kl_comparison}\subref{fig:complex_log_kl} plots the logarithm of the mean error, making the difference in late-time tracking accuracy more apparent. As summarized in Table~\ref{tab:crossing_tracking_metrics}, including \(\hat v_d\) gives a mean paired reduction of \(95.40\%\) in the late-time tracking error, with an improvement in all ten runs. These results show that accounting for the evolution of the task-demand density improves tracking beyond that achieved by the static-density feedback with a time-varying reference density.

\begin{table}[!htbp]
    \centering
    \caption{KL error statistics over ten random seeds. For each run, the late-time error is averaged over \(t\in[16.8,24]\) seconds, and the final error is evaluated at \(t=24\) seconds. Entries report the mean \(\pm\) sample standard deviation across the ten runs. Both numerical columns are in units of \(10^{-3}\).}
    \label{tab:crossing_tracking_metrics}
    \small
    \setlength{\tabcolsep}{4pt}
    \begin{tabular}{lcc}
        \hline
        Controller & Late-time mean & At \(t=24\)\\
        \hline
        With \(\hat v_d\) & \(9.101\pm0.190\) & \(9.434\pm0.273\)\\
        Without \(\hat v_d\) & \(197.777\pm1.030\) & \(153.814\pm1.155\)\\
        \hline
    \end{tabular}
\end{table}

\begin{figure*}[!t]
    \centering
    \captionsetup[subfigure]{skip=3pt}
    \begin{subfigure}{0.45\textwidth}
        \centering
        \includegraphics[width=\linewidth]{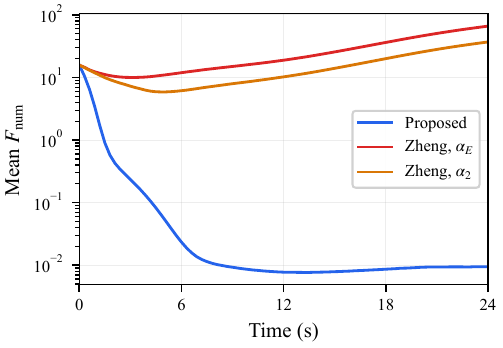}
        \caption{KL tracking error.}
        \label{fig:zheng_kl_panel}
    \end{subfigure}\hfill
    \begin{subfigure}{0.45\textwidth}
        \centering
        \includegraphics[width=\linewidth]{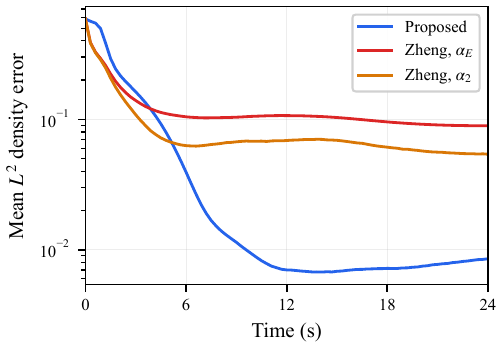}
        \caption{\(L^2\) density error.}
        \label{fig:zheng_l2_panel}
    \end{subfigure}
    \caption{Comparison with the adapted density-tracking feedback of Zheng~\cite{Zheng2021FieldEstimation} over evaluation seeds 23--32. The curves show pointwise means across the ten runs. Both vertical axes use logarithmic scales. The gains \(\alpha_E=3.81\) and \(\alpha_2=21.79\) are selected using calibration seeds 101--103 by energy bisection and local interval-halving refinement starting at \(\alpha=90\), respectively. All tested and selected gains use at most two decimal places.}
    \label{fig:zheng_tracking_comparison}

    \vspace{2mm}
\begin{minipage}{\textwidth}
    \centering
    \captionof{table}{Comparison over ten evaluation seeds. Late-time errors are means \(\pm\) sample standard deviations across runs; control energy and saturation fractions are averaged across runs. All controllers have the same speed limit.}
    \label{tab:zheng_tracking_metrics}
    \small
    \setlength{\tabcolsep}{4pt}
    \begin{tabular}{lccccc}
        \hline
        Controller & Gain & Late-time KL & Late-time \(L^2\) & \(\mathcal E\) & Saturation (\%)\\
        \hline
        Proposed & \(K=0.9\) & \(0.00910\pm0.00019\) & \(0.00767\pm0.00012\) & \(19.016\) & \(3.98\)\\
        Zheng's (energy) & \(\alpha=3.81\) & \(47.942\pm0.263\) & \(0.09230\pm0.00083\) & \(19.135\) & \(3.76\)\\
        Zheng's (refined \(L^2\)) & \(\alpha=21.79\) & \(26.608\pm0.277\) & \(0.05769\pm0.00072\) & \(168.575\) & \(73.72\)\\
        \hline
    \end{tabular}
\end{minipage}
\end{figure*}

\subsection{Performance Comparison }\label{subsec:simulation_zheng_comparison}

We further compare the proposed controller with the density-tracking feedback in Eq.~(8) of Zheng~\cite{Zheng2021FieldEstimation}. For the common diffusion coefficient \(D\), we adapt this feedback using the desired density \(p_r = \rho_d\) and the filtered transport estimate \(\hat v_d\). This supplies both controllers with the same
estimated target-motion information and gives
\begin{equation}
    \hat u_{\mathrm Z}={}-(\alpha-D)\nabla\ln\widetilde\rho_h^N+\frac{\rho_d}{\widetilde\rho_h^N} \bigl(\alpha\nabla\ln\rho_d+\hat v_d\bigr),
    \label{eq:zheng_comparison_control}
\end{equation}
where \(\alpha>0\) is a design parameter and \(\widetilde\rho_h^N\) is the full-space KDE in \eqref{eq:kde_estimate}. All controllers use the same target observations, KDE bandwidth, initial particles, Brownian increments, speed limit, and integration scheme.  In addition to the KL error, we evaluate the discrete \(L^2\) density error
\begin{equation}
    E_{2,\mathrm{num}}(t) :=\left[\sum_{p,q} \bigl(\bar\rho_{pq}(t)-\bar\rho_{d,pq}(t)\bigr)^2 \Delta A\right]^{1/2}.\label{eq:discrete_l2_metric}
\end{equation}

To compare tracking accuracy and control cost, we use the mean late-time \(L^2\) error and the applied control energy per agent, \(\mathcal E=(\Delta t/N)\sum_{k=0}^{T/\Delta t-1}\sum_{i=1}^N|u_{\mathrm a}(x_i^k,t_k)|^2\). We also report the fraction of agent-time samples at which the unsaturated drift exceeds \(u_{\max}\). 
Keeping \(K=0.9\), we choose both gains by bisection using calibration seeds 101--103, with all trial gains restricted to two decimal places.
\begin{enumerate}
    \item For comparable control cost, we select a gain whose mean applied energy approximately matches that of the proposed controller, giving \(\alpha=3.81\).

    \item For comparable tracking accuracy, we select the gain with the smallest mean late-time \(L^2\) error among the gains tested in this search, giving \(\alpha=21.79\).
\end{enumerate}
Both gains are then fixed for the evaluation seeds 23--32. Figure~\ref{fig:zheng_tracking_comparison} and Table~\ref{tab:zheng_tracking_metrics} show the results:
\begin{enumerate}
    \item At \(\alpha=3.81\), the mean applied control energies for Zheng's controller is \(19.135\), differing by \(0.63\%\) relative to the proposed controller value \(19.016\). At this similar energy level, the proposed controller achieves smaller late-time KL and \(L^2\) errors. The mean late-time \(L^2\) errors are \(0.00767\) and \(0.09230\), respectively, corresponding to a \(91.69\%\) reduction.

    \item At \(\alpha=21.79\), Zheng's controller has a mean late-time \(L^2\) error of \(0.05769\), which is \(7.52\) times the proposed value. Its mean applied energy is \(168.575\), approximately \(8.87\) times that of the proposed controller, and its drift is saturated for \(73.72\%\) of the agent-time samples.
\end{enumerate}
Thus, in this experiment, the proposed controller achieves smaller tracking errors at comparable applied energy, and even at the gain with the smallest mean late-time \(L^2\) error, Zheng’s controller has a larger error than the proposed controller.


\section{Conclusion}\label{section6}

This paper developed a mean-field feedback approach for cooperative tracking of multiple moving targets. We formulated multi-target tracking as a density-tracking problem by encoding the targets' tracking and sensing demands in a task-demand density. A feedback law was then derived from an Onsager-type variational problem, and the resulting closed-loop system was shown to admit a unique global weak solution with exponential KL convergence. For the estimated feedback law based on Gaussian KDE and filtered target velocities, we established a conditional KL error bound, thereby establishing the robustness of convergence with respect to the feedback mismatch. 

The simulations demonstrate the effectiveness and efficiency of our proposed approach: the swarm successfully tracks the moving targets, with the evolution of the mean KL error qualitatively consistent with the theoretical analysis, and the comparison with Zheng's controller shows smaller tracking errors at comparable applied energy. Our future work will investigate more general agent dynamics, including models with states evolving on \(\mathrm{SE}(n)\), and explore formulations of multi-target tracking based on information thermodynamics.

\section*{CRediT authorship contribution statement}

\textbf{Changming Lei:} Conceptualization, Methodology, Formal analysis, Software, Writing -- original draft.
\textbf{Junhao Dai:} Methodology, Validation, Writing -- review \& editing.
\textbf{Donghua Shi:} Methodology, Supervision, Writing -- review \& editing.

\bibliographystyle{elsarticle-num} 
\bibliography{ref.bib} 

\end{document}